\documentclass{amsart}
\usepackage{amssymb,amsmath,amsthm,amsbsy}
\usepackage[all]{xy}

\newtheorem{thm}{Theorem}[section]
\newtheorem{prop}[thm]{Proposition} 
\newtheorem{lemma}[thm]{Lemma}
\newtheorem{cor}[thm]{Corollary}

\theoremstyle{definition} 
\newtheorem{dfn}[thm]{Definition}

\theoremstyle{remark}
\newtheorem{rmk}[thm]{Remark}
\newtheorem{ex}[thm]{Example}

\newcommand{\comment}[1]{}

\newcommand{\la}[4]{
\xymatrix{#1 \ar[r] \ar@<2pt>[d] \ar@<-2pt>[d] & #2 \ar@<2pt>[d] \ar@<-2pt>[d] \\
#3 \ar[r] & #4}}
\newcommand{\dvb}[4]{
\xymatrix{#1 \ar[r] \ar[d] & #2 \ar[d] \\
#3 \ar[r] & #4}}

\newcommand{\pdiff}[2]{\frac{\partial #1}{\partial #2}}

\newcommand{\arrows}{\rightrightarrows}
\renewcommand{\to}{\longrightarrow}
\newcommand{\defequal}{:=}

\newcommand{\reals}{{\mathbb R}}

\newcommand{\integers}{{\mathbb Z}}

\newcommand{\lie}{\mathcal{L}}

\newcommand{\vect}{\mathfrak{X}}
\DeclareMathOperator{\sym}{S}
\newcommand{\bigS}{\mbox{{\large{$\sym$}}}}
\renewcommand{\bigwedge}{\mbox{{\Large{$\wedge$}}}}
\renewcommand{\tilde}[1]{\widetilde{#1}}
\renewcommand{\hat}[1]{\widehat{#1}}
\renewcommand{\bar}[1]{\overline{#1}}
\newcommand{\LA}{\mathcal{LA}}

\newcommand{\llbracket}{[\![}
\newcommand{\rrbracket}{]\!]}

\DeclareMathOperator{\im}{im}

\DeclareMathOperator{\ad}{ad}
\DeclareMathOperator{\Ad}{Ad}

\begin{document}
\title{$Q$-algebroids and their cohomology}
\author{Rajan Amit Mehta}
\begin{abstract}
A $Q$-algebroid is a Lie superalgebroid equipped with a compatible homological vector field and is the infinitesimal object corresponding to a $Q$-groupoid.  We associate to every $Q$-algebroid a double complex.  As a special case, we define the BRST model of a Lie algebroid, which generalizes the BRST model for equivariant cohomology.  We extend to this setting the Mathai-Quillen-Kalkman isomorphism of the BRST and Weil models, and we suggest a definition of a basic subcomplex which, however, requires a choice of a connection.  Other examples include Roytenberg's homological double of a Lie bialgebroid, Ginzburg's model of equivariant Lie algebroid cohomology, the double of a Lie algebroid matched pair, and $Q$-algebroids arising from lifted actions on Courant algebroids.
\end{abstract}

\maketitle
\tableofcontents

\section{Introduction}
Lie algebroids are known to appear in various guises.  In addition to the usual anchor-and-bracket definition, a Lie algebroid structure on a vector bundle $A \to M$ may be characterized as a linear Poisson structure on $A^*$, a Gerstenhaber algebra structure on $\bigwedge \Gamma(A)$, or a differential graded algebra structure on $\bigwedge \Gamma(A^*)$.  The latter two may be interpreted supergeometrically as a degree $-1$ Poisson structure on $[-1]A^*$ and a homological vector field (a \emph{$Q$-manifold} structure \cite{schwarz}) on $[-1]A$, respectively \cite{aksz, royt, vaintrob, voronov}.  

Since a bundle map $A \to A'$ does not, in general, induce a map of sections $\Gamma(A) \to \Gamma(A')$ or a dual bundle map $(A')^* \to A^*$, the differential graded algebra formulation has the distinction of being the point of view where the notion of morphism is obvious.  This fact has been useful in making sense of categorical constructs involving Lie algebroids, even though a Lie algebroid structure cannot be described abstractly in terms of objects and morphisms.  For example, Roytenberg \cite{royt} generalized the homological approach to Lie bialgebras \cite{ks1, lr} to give a homological compatibility condition for Lie bialgebroids.  More recently, Voronov \cite{voronov:mack} generalized this idea further, giving a homological compatibility condition for double Lie algebroids \cite{mac:dbl2}.

This paper deals with a closely related notion of categorical double, namely that of a \emph{$Q$-algebroid}.  A $Q$-algebroid may be thought of as a Lie algebroid in the category of $Q$-manifolds (or vice versa).  As mentioned above, certain $Q$-manifolds are associated to Lie algebroids.  $Q$-manifolds also arise from $L_\infty$-algebras \cite{lada-stasheff} and in relation to the quantization of gauge systems (see, e.g.\ \cite{henneaux-teitelboim}).

To any $Q$-algebroid, we can associate a double complex that is essentially the complex for Lie algebroid cohomology, ``twisted'' by the homological vector field.  A first example is the ``odd tangent algebroid'' $[-1]TA \to [-1]TM$, where $A \to M$ is a Lie algebroid.  In this case, the algebra of cochains is $\Omega([-1]A)$, the algebra of differential forms on the supermanifold $[-1]A$.  Since $[-1]A$ is a $Q$-manifold with homological vector field $d_A$ (the Lie algebroid differential), the Lie derivative with respect to $d_A$ is a differential operator on $\Omega([-1]A)$ that commutes with the de Rham operator; thus $\left(\Omega([-1]A), d, \lie_{d_A}\right)$ is a double complex.  When $A = M \times \mathfrak{g}$ is an action algebroid, this double complex is naturally isomorphic to the complex of the BRST model for equivariant cohomology \cite{kalkman:brst, osb}, which prompts us to view the double complex of $[-1]TA$ as a generalized BRST model.  Similarly, the complex $\left(\Omega([-1]A), d\right)$ may be viewed as a generalization of the Weil model\footnote{We are aware of independent work by Abad and Crainic \cite{abad-crainic}, where a Lie algebroid generalization of the Weil model is given without using the language of supergeometry.}.  In this setting, the isomorphism of Mathai-Quillen \cite{mq} and Kalkman \cite{kalkman:brst}, which relates the two models, takes a surprisingly simple, coordinate-free form.  Moreover, Kalkman's one-parameter family  \cite{kalkman:brst} that interpolates between the Weil and BRST models is immediately apparent in this context.

The Weil and BRST complexes are both models for $\Omega(M \times EG)$.  In order to obtain equivariant cohomology\footnote{As is usual for infinitesimal models of equivariant cohomology, we assume that the Lie group $G$ is compact and connected.}, one must first restrict to a subcomplex of forms that are basic for the $G$-action; it is this step that presents a problem in the general situation, as there does not seem to be a natural choice of a basic subcomplex of $\Omega([-1]A)$.  However, if we choose a linear connection on $A$, then we may define a basic subcomplex that, in the case of the canonical connection on $M \times \mathfrak{g}$, agrees with the usual subcomplex, which is just the Cartan model.  It remains unclear whether or not the basic cohomology depends in general on the choice of connection.

The notion of basic subcomplex extends to the double complex of a double Lie algebroid \cite{voronov} with \emph{generalized connection} \cite{alfonso:dvb}.  In the case of a vacant double Lie algebroid \cite{mac:dbl2} associated to a matched pair \cite{mokri} $(A,B)$ of Lie algebroids, there is a unique generalized connection, and the basic subcomplex computes the $A$-invariant Lie algebroid cohomology of $B$.  Examples of vacant double Lie algebroids include the Drinfel'd \cite{drinfeld:ybe} double $\mathfrak{g} \oplus \mathfrak{g}^*$ of a Lie bialgebra and the Lie algebroid $(M \times \mathfrak{g}) \oplus T^*M$ associated to a Poisson action \cite{lu}.  The latter arises in Mackenzie's \cite{mac:reduction} framework for Poisson reduction.

Other examples that fit into the framework of $Q$-algebroids include Roytenberg's \cite{royt} homological double of a Lie bialgebroid and Ginzburg's \cite{ginzburg} equivariant Lie algebroid cohomology.  As an example of a $Q$-algebroid that does not arise from a double Lie algebroid, we construct a $Q$-algebroid associated to a ``lifted action'' \cite{bcg} of a Lie algebra on a Courant algebroid.

$Q$-algebroids are the infinitesimal objects associated to $Q$-groupoids \cite{m:q}.  It was shown in \cite{m:q} that a double complex may be associated to a $Q$-groupoid by using the homological vector field to twist the smooth groupoid cohomology. It is natural to ask, then, what the relationship is between the cohomology of a $Q$-groupoid and that of its $Q$-algebroid.  As the above discussion of equivariant cohomology illustrates, the two are not generally isomorphic, but it may be possible to remedy the discrepancy by passing to a subcomplex at the $Q$-algebroid level.  In \cite{mw}, the relationship is explored further by showing that the van Est homomorphism \cite{crainic:vanest, wx} extends to a morphism of double complexes. 

The structure of the paper is as follows.  In \S\ref{sec:super} and \S\ref{sec:vfcart}, we provide a brief introduction to vector bundles and the Cartan calculus of differential forms on ($\integers$-graded) supermanifolds.  In \S\ref{section:superalg}, we define superalgebroids and show that the $[-1]$ functor from the category of Lie algebroids to the category of $Q$-manifolds extends to the category of superalgebroids.  As an application, we show in \S\ref{sec:poisson} that degree a $k$ Poisson structure on a supermanifold $\mathcal{M}$ may be associated to a superalgebroid structure on $[k]T^*\mathcal{M}$.  Section \ref{sec:morphic} deals with morphic vector fields and culminates in the definition of $Q$-algebroids.  The heart of the paper is \S\ref{sec:brst}, in which the double complex of a $Q$-algebroid is introduced and an in-depth study of $Q$-algebroids of the form $[-1]TA$ is undertaken.  Finally in \S\ref{sec:other}, we describe the other examples, mentioned above, of $Q$-algebroids and their double complexes.

Throughout the paper, we use calligraphic letters such as $\mathcal{M}$ to denote supergeometric objects and normal letters to denote objects that are assumed to be ordinary (i.e.\ not super).

\section{Vector bundles in the category of supermanifolds}\label{sec:super}
In this section, we describe some of the basic properties of vector bundles in the category of $\integers$-graded supermanifolds.  We will follow the notation of \cite{m:q} (also see \cite{mythesis}) regarding supermanifolds.

\subsection{Vector bundles}
Let $\mathcal{M}$ be a supermanifold with support $M$, and let $U$ be an open subset of $M$.  Denote by $\mathcal{M}|_U$ the supermanifold with support $U$ whose function sheaf is the restriction of $C^\infty(\mathcal{M})$ to $U$.

\begin{dfn} \label{dfn:bundle-top}A \emph{vector bundle} of rank $\{k_i\}$ over $\mathcal{M}$ is a supermanifold $\mathcal{E}$ and a surjection $\pi: \mathcal{E} \to \mathcal{M}$ equipped with an atlas of local trivializations $\mathcal{E}|_{\pi_0^{-1}(U)} \cong \mathcal{M}|_U \times \reals^{\{k_i\}}$ such that the transition map between any two local trivializations is linear in the fibre coordinates.
\end{dfn}
In particular, we note that a vector bundle has a well-defined subspace of linear functions $C^\infty_{lin}(\mathcal{E})$.  We will assume all vector bundles to have finite total rank, in that $\sum k_i < \infty$.

\begin{ex}[Trivial degree $j$ line bundles] Let $[j]\reals$ denote the coordinate superspace with one dimension in degree $j$.  Then $[j]\reals_\mathcal{M} \defequal \mathcal{M} \times [j]\reals$ is the trivial degree $j$ line bundle over $\mathcal{M}$.
\end{ex}

\subsection{Degree-shift functors}
For any integer $j$, there is a \emph{degree-shift functor} $[j]$ on the category of vector bundles, defined by
\begin{equation*}
[j]\mathcal{E} \defequal [j]\reals_\mathcal{M} \otimes \mathcal{E}.
\end{equation*}

\begin{rmk}
\begin{enumerate}
\item It is useful to consider the functor $[j]$ in terms of local fibre coordinates, as follows.  Let $\{\xi^a\}$ be local fibre coordinates on $\mathcal{E}$.  Then local fibre coordinates on $[j]\mathcal{E}$ are of the form $\{\xi^a \epsilon_j\}$, where $\epsilon_j$ is the fibre coordinate on $[j]\reals_\mathcal{M}$.  Note that $|\epsilon_j| = -j$.

\item It is clear from the definition that $[j][k]\mathcal{E}$ is isomorphic to $[k][j]\mathcal{E}$, but according to the sign rule, the canonical isomorphism is the one that sends $\xi^a \epsilon_k \epsilon_j$ to $(-1)^{jk} \xi^a \epsilon_j \epsilon_k$.  Similarly, $[j][k]\mathcal{E}$ is isomorphic to $[j+k]\mathcal{E}$, but the isomorphism requires a choice of identification of $[j]\reals \otimes [k]\reals$ with $[j+k]\reals$, which is only canonical up to a factor of $(-1)^{jk}$.  However, we may assume that such a choice has been fixed for all $j$ and $k$ (e.g.\ for $j<k$, identify $\epsilon_k \epsilon_j$ with $\epsilon_{j+k}$).
\end{enumerate}
\end{rmk}

\subsection{Sections}
Let $\pi: \mathcal{E} \to \mathcal{M}$ be a vector bundle. 

\begin{dfn} A degree $j$ \emph{section} of $\mathcal{E}$ is a map $X: \mathcal{M} \to [-j]\mathcal{E}$ such that $[-j]\pi \circ X = id_\mathcal{M}$.  The space of degree $j$ sections is denoted $\Gamma_j (\mathcal{E})$.  The \emph{space of sections} is $\Gamma(\mathcal{E}) \defequal \bigoplus_{j \in \integers} \Gamma_j (\mathcal{E})$.
\end{dfn}

\begin{rmk}
It is immediate from the definition that $\Gamma_{j-k}(\mathcal{E}) \cong \Gamma_j([k]\mathcal{E})$.  In other words, $\Gamma([k]\mathcal{E}) \cong [k]\Gamma(\mathcal{E})$, where the $[k]$ on the right-hand side is the usual degree-shift functor for graded vector spaces.  This is our main justification for why our notation for the degree-shift functor differs by a sign from, e.g.\ \cite{kont:deformation, royt:graded, severa}.
\end{rmk}

It is clear that, at least locally, a degree $j$ section $X$ is completely determined by the induced map $X^*: C^\infty_{lin}([-j]\mathcal{E}) \to C^\infty(\mathcal{M})$.
In a local trivialization, let $\{\xi^a\}$ be a set of fibre coordinates on $\mathcal{E}$, and let $\{\xi^a \epsilon_{-j}\}$ be the corresponding set of fibre coordinates on $[-j]\mathcal{E}$.  Then a set of local functions $\{f^a\} \in C^\infty(\mathcal{M})$ such that $|f^a| = |\xi^a| + j$ determines a local section $X \in \Gamma_j (\mathcal{E})$ with the property $X^*(\xi^a \epsilon_{-j}) = f^a$.  In particular, the \emph{frame of sections} $\{X_a\}$ dual to the coordinates $\{\xi^a\}$ is defined by the properties $|X_a| = -|\xi^a|$ and $X_a^*(\xi^b \epsilon_{|\xi^a|}) = \delta_a^b$.

Locally, any linear function $\alpha \in C^\infty_{lin}(\mathcal{E})$ can be written in the form $\alpha = \xi^a g_a$, where $g_a \in C^\infty(\mathcal{M})$.  The map $\xi^a g_a \mapsto \xi^a \epsilon_{-j} g_a$ gives a degree $j$, right $C^\infty(\mathcal{M})$-module isomorphism from $C^\infty_{lin}(\mathcal{E})$ to $C^\infty_{lin}([-j]\mathcal{E})$.  A left $C^\infty(\mathcal{M})$-module isomorphism may be defined similarly.  If $X \in \Gamma_j (\mathcal{E})$, then we can use the right module isomorphism
 of $C^\infty_{lin}(\mathcal{E})$ and $C^\infty_{lin}([-j]\mathcal{E})$ to obtain a degree $j$ map $X_0^*: C^\infty_{lin}(\mathcal{E}) \to C^\infty(\mathcal{M})$.  This map is a right $C^\infty(\mathcal{M})$-module homomorphism, and it satisfies the equation
\begin{equation*} \label{eqn:modhom}X_0^* (f\alpha) = (-1)^{|f|j} f X_0^*(\alpha) \end{equation*}
for any $f \in C^\infty(\mathcal{M})$ and $\alpha \in C^\infty_{lin}(\mathcal{E})$.

\begin{dfn} The pairing $\langle \cdot , \cdot \rangle: \Gamma(\mathcal{E}) \otimes C^\infty_{lin}(\mathcal{E}) \to C^\infty(\mathcal{M})$ is defined by 
\begin{equation}\label{eqn:pairing}\langle X, \alpha \rangle \defequal X_0^*(\alpha).\end{equation}
\end{dfn}

\begin{rmk}
It follows from the properties of $X_0^*$ that the pairing $\langle \cdot , \cdot \rangle$ is a right $C^\infty(\mathcal{M})$-module homomorphism, and one can see in local coordinates that the pairing is nondegenerate.  We may thus define a $C^\infty(\mathcal{M})$-module structure on $\Gamma(\mathcal{E})$ by the equation
\begin{equation*}\langle Xf, \alpha \rangle = \langle X, f \alpha \rangle. \end{equation*}
With this module structure, the pairing is a bimodule homomorphism.  
\end{rmk}

\subsection{Duals}\label{sub:duals}

Given a vector bundle $\mathcal{E} \to \mathcal{M}$, one can form a \emph{dual bundle} $\mathcal{E}^*$, where a local trivialization of $\mathcal{E}$ with fibre coordinates $\{\xi^a\}$ is associated to a local trivialization of $\mathcal{E}^*$ whose fibre coordinates are the dual frame of sections $\{X_a\}$, and transition maps for $\mathcal{E}$ determine transition maps for $\mathcal{E}^*$ such that the pairing (\ref{eqn:pairing}) is preserved.  By construction, we have that $\Gamma(\mathcal{E}^*) = C^\infty_{lin}(\mathcal{E})$.

The relationship between the grading of a vector bundle and its dual is opposite, in that if $\mathcal{E}$ is of rank $\{k_i\}$, then $\mathcal{E}^*$ is of rank $\{k_{-i}\}$.  In particular, $([j]\mathcal{E})^* = [-j]\mathcal{E}^*$.

\section{Vector fields and the Cartan calculus on supermanifolds}\label{sec:vfcart}

\subsection{Vector fields}\label{sub:vfields}

Let $\mathcal{M}$ be a supermanifold.

\begin{dfn} A \emph{vector field} of degree $j$ on $\mathcal{M}$ is a degree $j$ (left) derivation $\phi$ of $C^\infty(\mathcal{M})$, i.e.\ a linear operator such that, for any homogeneous functions $f, g \in C^\infty(\mathcal{M})$,
\begin{equation*} |\phi (f)| = j + |f| \end{equation*}
and 
\begin{equation*} \phi(fg) = \phi(f)g + (-1)^{j|f|}f \phi(g). \end{equation*}
\end{dfn}

The space of vector fields on $\mathcal{M}$ is denoted $\vect(\mathcal{M})$.  There is a natural $C^\infty(\mathcal{M})$-module structure on $\vect(\mathcal{M})$, defined by the property
\begin{equation*}
(f\phi)(g) = f \phi(g).
\end{equation*}

\begin{dfn}The \emph{Lie bracket} of two vector fields $\phi$ and $\psi$ is 
\begin{equation*}[\phi,\psi] = \phi\psi - (-1)^{|\phi||\psi|}\psi \phi.\end{equation*}
\end{dfn}

If $|\phi| = p$ and $|\psi|= q$, then $[\phi,\psi]$ is a vector field of degree $p+q$.  The Lie bracket gives $\vect(\mathcal{M})$ the structure of a Lie superalgebra; namely, for any homogeneous vector fields $\phi$, $\psi$, $\chi$ and function $f$:
\begin{enumerate}
\item $[\phi,\psi] = (-1)^{1+ |\phi||\psi|}[\psi,\phi]$ (antisymmetry),
\item $[\phi,f\psi]= \phi(f)\psi + (-1)^{|\phi||f|}f[\phi,\psi]$ (Leibniz rule),
\item $(-1)^{|\phi||\chi|}[\phi,[\psi,\chi]] + (-1)^{|\psi||\phi|}[\psi,[\chi,\phi]] + (-1)^{|\chi||\psi|}[\chi,[\phi,\psi]] = 0$ (Jacobi identity).
\end{enumerate}

In \S\ref{sec:diff}, the tangent bundle will be defined, and in \S\ref{sec:cartan} we will show that the space of vector fields can be identified with the space of sections of the tangent bundle.

\subsection{(Pseudo)differential forms}\label{sec:diff}

For a superdomain $\mathcal{U}$ of dimension $\{p_i\}$, the tangent bundle $T\mathcal{U}$ is the trivial bundle $\mathcal{U}\times\reals^{\{p_i\}}$.  Let $\{x^i\}$ be coordinates on $\mathcal{U}$, and denote the fibre coordinates by $\{\dot{x}^i\}$ (so that $|\dot{x}^i| = |x^i|$).  Naturally associated to any morphism $\mu: \mathcal{U} \to \mathcal{V}$ is a bundle map $T\mu$, defined as follows.  Let $\{y^i, \dot{y}^i\}$ be coordinates on $T\mathcal{V}$.  Then
\begin{equation} \label{eqn:tangent}
(T\mu)^* (\dot{y}^i) = \dot{x}^j \pdiff{}{x^j} [\mu^*(y^i)] .
\end{equation}
As supermanifolds are locally modeled on superdomains, these properties completely define the tangent functor on the category of supermanifolds.

Applying the $[-1]$ functor, we obtain the \emph{odd tangent bundle} $[-1]T\mathcal{M}$, where the fibre coordinates $\{\dot{x}^i\}$ are now considered to be of degree $|x^i| + 1$.  For an ordinary manifold $M$, it is clear from (\ref{eqn:tangent}) that the sheaf of functions on $[-1]TM$ is equal to $\Omega(M)$.  However, in the general case (specifically, when $\mathcal{M}$ has coordinates of degree $-1$), $[-1]T\mathcal{M}$ can have fibre coordinates of degree $0$ in which functions are not necessarily polynomial.  For this reason, $C^\infty([-1]T\mathcal{M})$ is called the algebra of \emph{pseudodifferential forms} \cite{bl2}.

\subsection{Cartan calculus on supermanifolds}\label{sec:cartan}

There is a canonical degree $1$ vector field $d$, known as the \emph{de Rham vector field}, on $[-1]T\mathcal{M}$ that is given locally by the formula
\begin{equation*} d = \dot{x}^i \pdiff{}{x^i}.\end{equation*}
It follows from (\ref{eqn:tangent}) that $d$ does not depend on the choice of coordinates.

For any $X \in \Gamma(T\mathcal{M})$, the identification $\Gamma([-1]T\mathcal{M})$ with $[-1]\Gamma(T\mathcal{M})$ allows us to identify $X$ with a degree $|X| - 1$ section $X_{-1} \in \Gamma([-1]T\mathcal{M})$.  

\begin{dfn}   The \emph{contraction} $\iota_X$ is the unique degree $|X| - 1$ vector field on $[-1]T\mathcal{M}$ satisfying, for any $\alpha \in C^\infty_{lin}([-1]T\mathcal{M})$,
\begin{equation}\label{eqn:contraction}\iota_X(\alpha) = \langle X_{-1}, \alpha \rangle. \end{equation}
\end{dfn}
Note that, if we view $f \in C^\infty(\mathcal{M})$ as a fibrewise-constant function on $[-1]T\mathcal{M}$, we have $\iota_X(f) = 0$.

\begin{rmk} \label{rmk:contraction} More generally, for any vector bundle $\mathcal{E}$, a section $X \in \Gamma(\mathcal{E})$ induces a degree $|X| - 1$ vector field $\iota_X \in \vect([-1]\mathcal{E})$, defined as in (\ref{eqn:contraction}).
\end{rmk}

\begin{dfn} Let $X \in \Gamma(T\mathcal{M})$.  The \emph{Lie derivative} with respect to $X$ is the degree $|X|$ vector field $\lie_X$ on $[-1]T\mathcal{M}$ defined by 
\begin{equation}\label{eqn:lie}\lie_X = [\iota_X, d].\end{equation}
\end{dfn}

\begin{prop}\label{prop:tangent}
The space of sections of the tangent bundle $T\mathcal{M}$ is naturally isomorphic to the space of vector fields on $\mathcal{M}$.
\end{prop}

\begin{proof} Let $X \in \Gamma(T\mathcal{M})$.  For any $f \in C^\infty(\mathcal{M})$, $\lie_X f = \left\langle X_{-1}, df\right\rangle \in C^\infty(\mathcal{M})$, so $\lie_X$ may be restricted to $C^\infty(\mathcal{M})$ to get a vector field on $\mathcal{M}$.

On the other hand, let $\phi \in \vect(\mathcal{M})$.  Then $\phi^V \in \vect([-1]T\mathcal{M})$, defined by the properties
\begin{align*}\phi^V (f) = 0, && \phi^V (df) = \phi(f) \end{align*}
for all $f \in C^\infty(\mathcal{M})$, restricts to a right $C^\infty(\mathcal{M})$-module morphism $C^\infty_{lin}([-1]T\mathcal{M}) \to C^\infty(\mathcal{M})$.  By the nondegeneracy of the pairing (\ref{eqn:pairing}), there is a section $X \in \Gamma(T\mathcal{M})$ such that $\phi^V = \iota_X$.

It is easy to check that the two maps are inverses of each other.
\end{proof}

\begin{ex}\label{ex:oddcot} Let $M$ be an ordinary manifold, and consider the \emph{odd cotangent bundle} $[-1]T^*M$.  Then $C^\infty([-1]T^*M) = \bigS([1]\vect(M))$ is the algebra of \emph{multivector fields} on $M$ and is denoted by $\vect^\bullet (M)$.  The Lie bracket operation on vector fields naturally extends to a graded biderivation $[\cdot, \cdot]: \vect^q(M) \otimes \vect^{q'}(M) \to \vect^{q + q' -1} (M)$, known as the \emph{Schouten bracket}. From the supergeometric perspective, we may view the Schouten bracket as the degree -1 Poisson bracket associated to the canonical degree $1$ symplectic structure on $[-1]T^*M$ (see \cite{ksm}, \cite{lz}).

If $\mathcal{M}$ is a supermanifold, then $C^\infty([-1]T^*\mathcal{M})$ is the algebra of \emph{pseudomultivector fields}.  Again, the Lie bracket operation on vector fields naturally extends as a biderivation to a bracket on the algebra of pseudomultivector fields, and this bracket may be viewed as the Poisson bracket associated to the canonical degree $1$ symplectic structure on $[-1]T^*\mathcal{M}$.
\end{ex}

The proof of the following proposition is left as an exercise for the reader.
\begin{prop}\label{prop:cartancomm}
For any vector fields $X, Y$ on $\mathcal{M}$, the following relations are satisfied:
\begin{align*}
[d,d] &=0, & [\iota_X, \iota_Y] &= 0, \\
[\lie_X, \iota_Y] &= \iota_{[X,Y]}, & [\lie_X, \lie_Y] &= \lie_{[X,Y]}, \\
[d, \lie_X] &= 0.
\end{align*}
\end{prop}

From the above commutation relations, a Cartan-type formula for the de Rham differential may be derived.  
In particular, from the equations
\begin{align*}
\iota_X d &= \lie_X - (-1)^{|X|} d \iota_X, \\
\iota_Y \lie_X &= (-1)^{|X|(|Y| - 1)}(\lie_X \iota_Y - \iota_{[X,Y]}),
\end{align*}
it follows that for any $p$-form $\omega$ and any vector fields $X_0, \dots, X_p$,
\begin{align}
\begin{split}\label{eqn:cartdiff}
& \iota_{X_p} \cdots \iota_{X_0} d \omega = \sum_{i=0}^p (-1)^{i+ \sum_{j<i}|X_j|} \iota_{X_p} \cdots \lie_{X_i} \cdots \iota_{X_0}\omega \\
&= \sum_{i=0}^p (-1)^{i+\sum_{k<i}|X_k| + |X_i|\left(\sum_{k > i} (|X_k| - 1)\right)} \lie_{X_i}\iota_{X_p} \cdots \widehat{\iota_{X_i}} \cdots \iota_{X_0} \omega \\
& \quad + \sum_{j>i} (-1)^{j+\sum_{k<i}|X_k| + |X_i|\left(\sum_{j \geq k > i}(|X_k| - 1)\right)} \iota_{X_p} \cdots \iota_{[X_i,X_j]} \cdots \widehat{\iota_{X_i}} \cdots \iota_{X_0} \omega.
\end{split}
\end{align}

\section{Superalgebroids and $Q$-algebroids}
\label{section:superalg}

\begin{dfn}
Let $\mathcal{M}$ be a supermanifold, and let $\mathcal{A}$ be a supervector bundle over $\mathcal{M}$.  A \emph{superalgebroid} structure on $\mathcal{A}$ is a degree $0$ bundle map $\rho: \Gamma(\mathcal{A}) \to \vect(\mathcal{M})$ (the \emph{anchor}) and a Lie superbracket $[\cdot,\cdot]: \Gamma (\mathcal{A})  \times \Gamma (\mathcal{A})  \to \Gamma (\mathcal{A}) $ such that
\begin{equation*}
[X, fY] = \lie_{\rho(X)}f \cdot Y + (-1)^{|X||f|}f[X,Y]
\end{equation*}
for all $X, Y \in \Gamma (\mathcal{A})$ and $f \in C^\infty(\mathcal{M})$.
\end{dfn}

\begin{rmk} As in the ordinary case, the \emph{anchor identity} $\rho([X,Y]) = [\rho(X),\rho(Y)]$ may be deduced from the properties of a superalgebroid.
\end{rmk}

Many common constructions of Lie algebroids extend in a straightforward manner to the supergeometric situation.  In particular, examples of superalgebroids include the tangent bundle $T\mathcal{M}$ of a supermanifold $\mathcal{M}$, the action algebroid $\mathcal{M} \times \mathfrak{g} \to \mathcal{M}$ arising from an infinitesimal action of a Lie superalgebra $\mathfrak{g}$ on a supermanifold $\mathcal{M}$, and the cotangent bundle of a supermanifold with a degree $0$ Poisson structure.  A more interesting generalization occurs when one wishes to construct a superalgebroid associated to a Poisson structure of nonzero degree; we describe such a construction in \S\ref{sec:poisson}, after a discussion of differentials.

We will generally drop the prefix ``super-'' in the remainder of the paper, except for emphasis.  In the case of a superalgebroid, we will simply write ``algebroid'', reserving the modifier ``Lie'' for the ordinary case.

\subsection{Algebroids and $Q$-manifolds}\label{sub:algq}

As in the ordinary case, an algebroid structure on $\mathcal{A}$ corresponds to a differential $d_{\mathcal{A}} : \bigS^\bullet [1]\Gamma (\mathcal{A}^*) \to \bigS^{\bullet + 1} [1]\Gamma (\mathcal{A}^*)$, which may be viewed as a homological vector field on $[-1]\mathcal{A}$.  We will define $d_{\mathcal{A}}$ in a coordinate frame and then derive Cartan commutation relations, from which an invariant formula may be obtained.

Let $\{x^i\}$ be local coordinates on $\mathcal{M}$, and let $\{X_\alpha\}$ be a frame of sections of $\mathcal{A}$ where $|X_\alpha| = p_\alpha$.  Then the anchor may be described locally in terms of functions $\rho_\alpha^i$ such that $\rho (X_\alpha) = \rho_\alpha^i \pdiff{}{x^i}$.  Similarly, the bracket may be described in terms of the \emph{structure functions} $c_{\alpha \beta}^\gamma$, where $[X_\alpha, X_\beta] = c_{\alpha \beta}^\gamma X_\gamma$.  

If $\{\lambda^\alpha\}$ are the degree-shifted (in that $|\lambda^\alpha| = -p_\alpha + 1$) fibre coordinates dual to $\{X_\alpha\}$, $d_{\mathcal{A}}$ is defined locally as
\begin{equation}\label{eqn:dalocal}
d_{\mathcal{A}} = \lambda^\alpha \rho_\alpha^i \pdiff{}{x^i} - (-1)^{p_\alpha (p_\beta - 1)}\frac{1}{2} \lambda^\alpha \lambda^\beta c_{\alpha \beta}^\gamma \pdiff{}{\lambda^\gamma}.
\end{equation}

\begin{lemma}\label{lemma:dsquared}
$d_\mathcal{A}^2 = 0$.
\end{lemma}
\begin{proof}
A direct calculation shows that the vanishing of $d_\mathcal{A}^2 (x^i)$ for all $i$ is equivalent to the anchor identity, and the vanishing of $d_\mathcal{A}^2(\lambda^\alpha)$ for all $\alpha$ is equivalent to the Jacobi identity.
\end{proof}

For $X \in \Gamma(\mathcal{A})$, we have the contraction $\iota_X \in \vect([-1]\mathcal{A})$ (see Remark \ref{rmk:contraction}).  Additionally, we may define Lie derivatives as in (\ref{eqn:lie}) by $\lie_X \defequal [\iota_X, d_\mathcal{A}]$.  An immediate consequence of the Jacobi identity is that $[\lie_X, d_\mathcal{A}] = 0$.  The following identity is a nontrivial property of (\ref{eqn:dalocal}), but it can be easily verified in coordinates.

\begin{lemma}\label{lemma:liecontract}
For any $X,Y \in \Gamma(\mathcal{A})$, 
\begin{equation*}
[\lie_X, \iota_Y] = \iota_{[X,Y]}.
\end{equation*}
\end{lemma}
From Lemma \ref{lemma:liecontract}, it follows that $[\lie_X, \lie_Y] = \lie_{[X,Y]}$.  In summary, we have obtained commutation relations that are identical to those of Proposition \ref{prop:cartancomm} (with $d_\mathcal{A}$ in the place of $d$), and therefore the Cartan formula (\ref{eqn:cartdiff}) holds for $d_\mathcal{A}$.  Since $\lie_X$, when restricted to $C^\infty(\mathcal{M})$, acts as $\rho(X)$, we see that (\ref{eqn:cartdiff}) gives a formula for $d_\mathcal{A}$ that does not depend on a choice of coordinates.

In the other direction, $d_{\mathcal{A}}$ completely encodes the algebroid structure of $\mathcal{A}$ in the following sense.

\begin{thm}\label{thm:deriv}
Let $\mathcal{A} \to \mathcal{M}$ be a vector bundle.  Given a degree $1$ derivation $d_{\mathcal{A}}: \bigS^\bullet [1]\Gamma (\mathcal{A}^*) \to \bigS^{\bullet + 1} [1]\Gamma (\mathcal{A}^*)$, there is a unique bundle map $\rho: \Gamma(\mathcal{A}) \to \vect(\mathcal{M})$ and a unique graded skew-symmetric bilinear map $[\cdot,\cdot]: \Gamma (\mathcal{A}) \otimes \Gamma (\mathcal{A}) \to \Gamma (\mathcal{A}) $, satifying the Leibniz rule, such that (\ref{eqn:cartdiff}) holds.  The Jacobi identity is satisfied if and only if $d_{\mathcal{A}}^2 = 0$.
\end{thm}

\begin{proof}
Applying (\ref{eqn:cartdiff}) to $f \in C^\infty(\mathcal{M})$ and $\omega \in \Gamma (\mathcal{A}^*)$, we obtain the following equations, which uniquely define $\rho$ and $[\cdot,\cdot]$:
\begin{align}\label{eqn:anchor}
\rho(X)f &= \iota_X d_{\mathcal{A}} f = [\iota_X, d_{\mathcal{A}}] f, \\
\begin{split}\label{eqn:bracket}
\iota_{[X,Y]} \omega &=  \rho(X) \iota_Y \omega - (-1)^{|X||Y|} \rho(Y)\iota_X \omega - (-1)^{|X|(|Y|-1)}\iota_Y \iota_X d_{\mathcal{A}} \omega \\
&= \left[ [\iota_X, d_{\mathcal{A}}], \iota_Y \right] \omega.
\end{split}
\end{align}
The last part of the proposition follows from the proof of Lemma \ref{lemma:dsquared}.
\end{proof}

Generalizing the role that $[-1]T\mathcal{M}$ plays in the theory of differential forms, $C^\infty([-1]\mathcal{A})$ may be called the algebra of \emph{algebroid pseudoforms}.  The vector field $d_{\mathcal{A}}$ is a homological vector field on $[-1]\mathcal{A}$; in other words, $[-1]\mathcal{A}$ naturally possesses the structure of a \emph{$Q$-manifold}.

Morphisms of algebroids may be defined in terms of the differentials; specifically, if $\mathcal{A} \to \mathcal{M}$ and $\mathcal{A}' \to \mathcal{N}$ are algebroids, then a smooth linear map $\tau: \mathcal{A} \to \mathcal{A}'$ is an algebroid morphism if $d_{\mathcal{A}}$ and $d_{\mathcal{A}'}$ are $[-1]\tau$-related, or, in other words, if $[-1]\tau$ is a morphism in the category of $Q$-manifolds.

We have thus extended to the case of superalgebroids the following result of \cite{vaintrob} (also see \cite{aksz}):
\begin{thm}\label{thm:minusone}
There is a functor $[-1]$ from the category of algebroids to the category of $Q$-manifolds.
\end{thm}

\begin{rmk}\label{rmk:anti}The $Q$-manifolds in the image of the $[-1]$ functor take a specific form; they are vector bundles such that their homological vector fields are quadratic in the fibre coordinates.  Following Voronov \cite{voronov}, we refer to such objects as \emph{anti-algebroids}.  The $[-1]$ functor is an equivalence of categories from the category of algebroids to the category of anti-algebroids, for which the inverse functor is $[1]$.
\end{rmk}

\subsection{Degree $k$ Poisson structures}\label{sec:poisson}

As an application of the homological approach to superalgebroids, we describe how a degree $k$ Poisson structure on a supermanifold $\mathcal{M}$ induces an algebroid structure on $[k]T^*\mathcal{M}$.  

\begin{dfn} A \emph{degree $k$ Poisson structure} on $\mathcal{M}$ is a degree $k$ bilinear map $\{\cdot,\cdot\}: C^\infty(\mathcal{M}) \otimes C^\infty(\mathcal{M}) \to C^\infty(\mathcal{M})$ such that, for all $f,g,h \in C^\infty(\mathcal{M})$,
\begin{enumerate}
\item $\{f,g\} = (-1)^{1+(|f|-k)(|g|-k)}\{g,f\}$,
\item $\{f,gh\} = \{f,g\}h + (-1)^{(|f|-k)|g|}g\{f,h\}$,
\item $\{f, \{g,h\}\} = \{ \{f,g\},h\} + (-1)^{(|f|-k)(|g|-k)}\{g,\{f,h\}\}$.
\end{enumerate}
\end{dfn}

As in the ordinary case, it may be seen that properties (1) and (2) are equivalent to the existence of a ``bivector field'' $\pi$, which in general is a degree $2-k$ element of $\bigS^2([1-k]\vect(\mathcal{M}))$, generating the Poisson bracket via the derived bracket formula
\begin{equation*}\{f,g\} = \left[ [f,\pi]_k,g \right]_k. \end{equation*}
Here, $[\cdot,\cdot]_k$ is the degree $k-1$ Schouten bracket on $C^\infty([k-1]T^*\mathcal{M})$.  Property (3) is then equivalent to the ``integrability condition'' $[\pi,\pi]_k = 0$.

The integrability condition implies that the operator $d_\pi \defequal [\pi, \cdot]_k$ is a homological vector field on $[k-1]T^*\mathcal{M}$, so $d_\pi$ gives $[k-1]T^*\mathcal{M}$ the structure of an anti-algebroid.  The anchor $\pi^\sharp$ and bracket $[\cdot,\cdot]_\pi$ for the associated algebroid structure on $[k]T^*\mathcal{M}$ may then be obtained from (\ref{eqn:anchor}) and (\ref{eqn:bracket}).  In particular, if we denote by $d_k$ the ``degree $k$ de Rham operator'' that, in particular, sends $C^\infty(\mathcal{M})$ to $\Gamma([k]T^*\mathcal{M}) = [k]\Omega^1(\mathcal{M})$, then $\pi^\sharp$ and $[\cdot,\cdot]_\pi$ may be shown to satisfy the properties
\begin{align*}
\pi^\sharp (d_k f)(g) &= \{f,g\}, \\
[d_k f, d_k g]_\pi &= d_k \{f,g\},
\end{align*}
for all $f,g \in C^\infty(\mathcal{M})$.

\subsection{Morphic vector fields}\label{sec:morphic}

In this section, we define morphic vector fields and describe their basic properties.  Morphic vector fields correspond, in the ordinary case, to infinitesimal algebroid automorphisms \cite{mac-xu}.  

Let $\mathcal{E} \to \mathcal{M}$ be a vector bundle.  

\begin{dfn}\label{dfn:linvf}A vector field $\Xi$ on $\mathcal{E}$ is \emph{linear} if, for every $\alpha \in C^\infty_{lin}(\mathcal{E})$, $\Xi(\alpha) \in C^\infty_{lin}(\mathcal{E})$.  The space of linear vector fields is denoted $\vect_{lin}(\mathcal{E})$.
\end{dfn}

\begin{rmk}It follows from the derivation property that, if $\Xi$ is a linear vector field on $\mathcal{E}$, then $C^\infty(\mathcal{M})$ (considered as a subalgebra of $C^\infty(\mathcal{E})$) is invariant under the action of $\Xi$ as a derivation.  It follows that there exists a unique vector field $\phi$ on $\mathcal{M}$ that is $\pi$-related to $\Xi$.  We call $\phi$ the \emph{base vector field}.
\end{rmk} 

\begin{prop}\label{prop:liniso} For all $j$, $\vect_{lin}(\mathcal{E})$ and $\vect_{lin}([-j]\mathcal{E})$ are naturally isomorphic to each other as $C^\infty(\mathcal{M})$-modules.
\end{prop}

\begin{proof}
Let $\Xi \in \vect_{lin}(\mathcal{E})$.  Using the left $C^\infty(\mathcal{M})$-module isomorphism $C^\infty_{lin}(\mathcal{E}) \to C^\infty_{lin}([-j]\mathcal{E})$, $\alpha \mapsto \hat{\alpha}$, we may define a linear operator $\Xi_j$ on $C^\infty_{lin}([-j]\mathcal{E})$ by the property
\begin{equation}\label{eqn:hatmap}
\Xi_j (\hat{\alpha}) = \widehat{\Xi(\alpha)}.
\end{equation}
From (\ref{eqn:hatmap}), the left $C^\infty(\mathcal{M})$-module isomorphism property of $\alpha \mapsto \hat{\alpha}$, and the fact that $\Xi$ is a linear vector field, it is simple to see that, for $f \in C^\infty(\mathcal{M})$,
\begin{equation}\label{eqn:opdev}
\Xi_j (f \hat{\alpha}) = \phi (f) \hat{\alpha} + (-1)^{|\Xi||f|} f \Xi_j (\hat{\alpha}).
\end{equation}
Equation (\ref{eqn:opdev}) implies that $\Xi_j$ may be extended as a derivation to a linear vector field on $[-j]\mathcal{E}$ of the same degree, and with the same base vector field, as $\Xi$.

It is immediate from the definition that $(f\Xi)_j = f \Xi_j$.  Since the map $\Xi \mapsto \Xi_j$ is a degree $0$ map, the fact that it is a left module isomorphism implies that it is also a right module isomorphism.
\end{proof}

\begin{rmk} 
In local coordinates, where $\{x^i\}$ are coordinates on $\mathcal{M}$ and $\{\xi^a\}$ are fibre coordinates, we may write 
\[\Xi = \phi^i (x) \pdiff{}{x^i} + \varphi^a_b(x) \xi^b \pdiff{}{\xi^a}.\]
Then the local description of $\Xi_j$ is simply
\[\Xi_j = \phi^i (x) \pdiff{}{x^i} + \varphi^a_b(x) \hat{\xi}^b \pdiff{}{\hat{\xi}^a}.\]
\end{rmk}

\begin{prop}\label{prop:hathom}
The map $\Xi \mapsto \Xi_j$ is a Lie algebra homomorphism.
\end{prop}

\begin{proof}
Since $\Xi_j$ is defined by the property (\ref{eqn:hatmap}), we have that
\begin{equation*}
\begin{split}
[\Xi_j, \Xi'_j] (\hat{\alpha}) &= \Xi_j \Xi_j^\prime (\hat{\alpha}) - (-1)^{|\Xi||\Xi'|}\Xi_j^\prime \Xi_j (\hat{\alpha}) \\
&= \Xi_j \left(\widehat{\Xi^\prime (\alpha)}\right) - (-1)^{|\Xi||\Xi'|}\Xi^\prime_j \left(\widehat{\Xi (\alpha)}\right) \\
&= \widehat{\Xi \Xi^\prime (\alpha)} - (-1)^{|\Xi||\Xi'|}\widehat{\Xi^\prime \Xi (\alpha)} \\
&= \widehat{[\Xi, \Xi'] (\alpha)},
\end{split}
\end{equation*}
from which it follows that $[\Xi_j, \Xi'_j] = [\Xi,\Xi']_j$.
\end{proof}

\begin{dfn} Let $\mathcal{A} \to \mathcal{M}$ be an algebroid.  A linear vector field $\Xi$ on $\mathcal{A}$ is called \emph{morphic} if $[d_\mathcal{A}, \Xi_1] = 0$.
\end{dfn}

\begin{rmk}When there is no possibility of ambiguity, $\Xi_1 \in \vect([-1]\mathcal{A})$ will also be called morphic if $\Xi$ is morphic.
\end{rmk}

\begin{rmk} For any $X \in \Gamma(\mathcal{A})$, the vector field $\lie_X \in \vect([-1]\mathcal{A})$ is morphic.  When $\mathcal{A} = T\mathcal{M}$, every morphic vector field is of the form $\lie_X$ for some $X \in \vect(\mathcal{M})$.  However, an algebroid may in general have ``outer symmetries'' that do not arise from sections.  An extreme example is an algebroid whose bracket is always $0$.  For such an algebroid, $d_{\mathcal{A}}=0$, so any linear vector field is morphic.  However, in this case $\lie_X=0$ for all $X \in \Gamma(\mathcal{A})$.
\end{rmk}

\begin{dfn} Let $\Xi$ be a morphic vector field on an algebroid $\mathcal{A}$.  We define the operator $D_\Xi$ on $\Gamma(\mathcal{A})$ by the equation
\begin{equation}\label{eqn:morph}
\iota_{D_\Xi X} = [\Xi_1, \iota_X].
\end{equation}
\end{dfn}

From (\ref{eqn:morph}), the following formula may be derived:
\begin{equation}\label{eqn:liemorph}
\begin{split}
\iota_{X_q} \cdots \iota_{X_1} \Xi_1 \omega =& 
(-1)^{|\Xi|\sum_{k=1}^q (|X_k| - 1)} \phi(\iota_{X_q} \cdots \iota_{X_1} \omega)\\
 &- \sum_{i=1}^q (-1)^{|\Xi|\sum_{k=1}^i (|X_k| - 1)}\iota_{X_q} \cdots \iota_{X_{i+1}} \iota_{D_\Xi X_i} \iota_{X_{i-1}} \cdots \iota_{X_1} \omega
 \end{split}
 \end{equation}
The reader may compare (\ref{eqn:liemorph}) with a similar formula for Lie derivatives in standard differential calculus. 

\begin{prop} \label{prop:morphbrack} Let $\mathcal{A} \to \mathcal{M}$ be an algebroid.  Then
\begin{enumerate}
\item The space of morphic vector fields on $\mathcal{A}$ is closed under the Lie bracket.
\item If $\Xi$ is a morphic vector field and $X \in \Gamma(\mathcal{A})$ then $[\Xi_1, \lie_X] = \lie_{D_\Xi X}$.
\end{enumerate}
\end{prop}

\begin{proof}
By the Jacobi identity, it is clear that if $\Xi$ and $\Xi'$ are linear vector fields such that $[\Xi_1, d_{\mathcal{A}}] = [\Xi^\prime_1, d_{\mathcal{A}}] = 0$, then $\left[ [\Xi_1, \Xi^\prime_1], d_{\mathcal{A}}\right] = 0$.  Using Proposition \ref{prop:hathom}, we conclude that $[\Xi, \Xi']$ is a morphic vector field.

For the second statement, the Jacobi identity and the fact that $[\Xi_1, d_{\mathcal{A}}] = 0$ imply that 
\[ \left[\Xi_1, [\iota_X, d_{\mathcal{A}}]\right] = \left[[\Xi_1, \iota_X], d_{\mathcal{A}}\right], \]
or, simply, $[\Xi_1, \lie_X] = \lie_{D_\Xi X}$.
\end{proof}

\begin{prop}\label{prop:morphop}Let $\Xi$ be a morphic vector field with base vector field $\phi$.  Then 
\begin{enumerate}
\item For any $X \in \Gamma(\mathcal{A})$, $\rho(D_\Xi X) = [\phi, \rho(X)]$,
\item $D_\Xi$ is a derivation with respect to the Lie bracket.
\end{enumerate}
\end{prop}
\begin{proof}
The first statement follows from the second part of Proposition \ref{prop:morphbrack}, restricted to $C^\infty(\mathcal{M})$.

For the second statement, we see that
\begin{equation}\begin{split}
\iota_{D_\Xi [X,Y]} &= \left[\Xi_1, [\lie_X, \iota_Y]\right]\\
 &= \left[ [\Xi_1, \lie_X], \iota_Y \right] + (-1)^{|\Xi||X|} \left[\lie_X, [\Xi_1, \iota_Y] \right] \\
 &= \iota_{[D_\Xi X, Y]} + (-1)^{|\Xi||X|} \iota_{[X, D_\Xi Y]}.\qedhere
\end{split}\end{equation}
\end{proof}

\begin{lemma}\label{lem:morphderiv}
Let $D$ be a linear operator on $\Gamma (\mathcal{A})$ and let $\phi \in \vect(\mathcal{M})$ such that
\begin{enumerate} 
\item $D (fX) = \phi(f)X + (-1)^{|f||D|}f D(X)$ for all $f \in C^\infty(\mathcal{M})$ and $X \in \Gamma(\mathcal{A})$, and
\item $D$ is a derivation of the Lie bracket.
\end{enumerate}
Then, for any $X \in \Gamma(\mathcal{A})$, $\rho(D(X)) = [\phi, \rho(X)]$.
\end{lemma}
\begin{proof}Using the Leibniz rule, we have that for any $X, Y \in \Gamma(\mathcal{A})$ and $f \in C^\infty(\mathcal{M})$,
\begin{equation}\begin{split}\label{morph1}
D \left( [X, fY] \right) =& D\left(\rho(X)(f) \cdot Y + (-1)^{|X||f|} f[X,Y] \right) \\
=& \phi \circ \rho(X)(f) \cdot Y + (-1)^{|D|(|X| + |f|)} \rho(X)(f) \cdot D(Y) + (-1)^{|X||f|} \phi(f) [X,Y] \\
&+ (-1)^{(|D| + |X|)|f|} f[D(X), Y] + (-1)^{(|D| + |X|)|f| + |D||X|} f[X,D(Y)].
\end{split}\end{equation}
On the other hand, 
\begin{equation}\begin{split}\label{morph2}
D \left( [X, fY] \right) =& [D(X),fY] + (-1)^{|D||X|} [X,D(fY)] \\
=& [D(X),fY] + (-1)^{|D||X|}[X, \phi(f)\cdot Y + (-1)^{|D||f|} f D(Y)] \\
=& \rho(D(X))(f) \cdot Y + (-1)^{(|D| + |X|)|f|} f[D(X),Y] + (-1)^{|D||X|} \rho(X) \circ \phi(f) \cdot Y \\
&+ (-1)^{|X||f|} \phi(f)[X,Y] + (-1)^{|D|(|X| + |f|)}\rho(X)(f) \cdot D(Y) \\
&+ (-1)^{|D|(|X|+|f|) + |X||f|} f[X,D(Y)].
\end{split}\end{equation}
After equating the results of (\ref{morph1}) and (\ref{morph2}) and cancelling terms, we obtain the equation
\begin{equation*}
\phi \circ \rho(X)(f) \cdot Y = \rho(D(X))(f) \cdot Y + (-1)^{|D||X|} \rho(X) \circ \phi(f) \cdot Y,
\end{equation*}
which may be written more simply as
\begin{equation}\label{morph3}
[\phi, \rho(X)] (f) \cdot Y = \rho(D(X))(f) \cdot Y.
\end{equation} 
Since (\ref{morph3}) holds for all $f$ and $Y$, we conclude that $[\phi, \rho(X)] = \rho(D(X))$.
\end{proof}

\begin{thm}\label{thm:morphderiv}
Let $D$ be a linear operator on $\Gamma (\mathcal{A})$ satisfying the hypotheses of Lemma \ref{lem:morphderiv}.  Then there exists a morphic vector field $\Xi \in \vect(\mathcal{A})$ such that $D = D_\Xi$.
\end{thm}
\begin{proof} Given such a $D$, let $\Xi_1$ be the degree $|D|$ linear operator on $C^\infty_{lin}([-1]\mathcal{A})$ defined by the property
\begin{equation*} \iota_X \Xi_1(\alpha) = (-1)^{|D|(|X|-1)} \left(\phi(\iota_X \alpha) - \iota_{D(X)} \alpha\right),
\end{equation*}
for all $X \in \Gamma(\mathcal{A})$ and $\alpha \in C^\infty_{lin}([-1]\mathcal{A})$.  Since, for any $f \in C^\infty(\mathcal{M})$, 
\begin{equation*}
\iota_X \Xi_1(f\alpha) = \iota_X \left( \phi(f) \alpha + (-1)^{|D||f|} f \Xi_1(\alpha) \right),
\end{equation*}
we may extend $\Xi_1$ uniquely to a degree $|D|$ linear vector field on $[-1]\mathcal{A}$ with base vector field $\phi$.

To show that $\Xi$ is morphic, it is sufficient to check that $[d_{\mathcal{A}}, \Xi_1] (\alpha) = 0$ for all $\alpha \in C^\infty_{lin}([-1]\mathcal{A})$.  A somewhat long, but direct, calculation using (\ref{eqn:cartdiff}) and (\ref{eqn:liemorph}) reveals that, for any $X, Y \in \Gamma(\mathcal{A})$,
\begin{equation*}\begin{split}
\iota_Y \iota_X [d_{\mathcal{A}}, \Xi_1] (\alpha) = &\pm \left( \rho(D(X)) - \phi \circ \rho(X) + (-1)^{|D||X|}\rho(X) \circ \phi \right) \iota_Y \alpha \\
& \pm  \left( \rho(D(Y)) - \phi \circ \rho(Y) + (-1)^{|D||Y|}\rho(Y) \circ \phi \right) \iota_X \alpha \\
& \pm \left( \iota_{D[X,Y]} - \iota_{[D(X),Y]} - (-1)^{|X||D|} \iota_{[X, D(Y)]} \right) \alpha,
\end{split}
\end{equation*}
where the signs for each line are omitted.  The first two lines vanish by Lemma \ref{lem:morphderiv}, and the final line vanishes because $D$ is by assumption a derivation of the Lie bracket.  It follows that $[d_{\mathcal{A}}, \Xi_1] = 0$.
\end{proof}

\begin{dfn} An algebroid equipped with a homological morphic vector field is called a \emph{Q-algebroid}.
\end{dfn}

\section{The BRST complex}\label{sec:brst}

\subsection{The cohomology of $Q$-algebroids}

If $(\mathcal{A} \to \mathcal{M}, \Xi)$ is a $Q$-algebroid, then $\Xi_1$ and $d_\mathcal{A}$ form a pair of commuting homological vector fields on $[-1]\mathcal{A}$.  Furthermore, the grading on the polynomial functions $\bigS [1]\Gamma(\mathcal{A}^*) \subseteq C^\infty([-1]\mathcal{A})$ splits into a double grading $(p,q)$, where $p$ is the ``cohomological grading'' and $q$ is the ``supermanifold grading''.  Specifically, a simple element of $\bigS [1]\Gamma(\mathcal{A}^*)$ is of the form $\hat{\alpha}_1\cdots\hat{\alpha}_p$, where $\alpha_i \in \Gamma(\mathcal{A}^*)$, and the hat denotes that the grading has been increased by $1$.  For such an element, the cohomological grading is $p$, and the supermanifold grading is $q \defequal \sum |\alpha_i|$.  With respect to the double grading, $\Xi_1$ is of degree $(0,1)$ and $d_\mathcal{A}$ is of degree $(1,0)$.  Therefore we have a double complex $\left(C^{p,q}(\mathcal{A}), \Xi_1, d_\mathcal{A}\right)$, where $C^{p,q}(\mathcal{A})$ denotes the space of degree $(p,q)$ elements of $\bigS [1]\Gamma(\mathcal{A}^*)$.

\begin{dfn}
The \emph{$Q$-algebroid cohomology} of $\mathcal{A}$ is the total cohomology of the associated double complex $\left(C^{p,q}(\mathcal{A}), \Xi_1, d_\mathcal{A}\right)$.
\end{dfn}

\begin{rmk}
Voronov \cite{voronov:mack} has recently shown that the structure of a double Lie algebroid \cite{mac:dbl2, mac:dbl}
\begin{equation}\label{eqn:dvb}
\dvb{D}{A}{B}{M}
\end{equation}
can be encoded in a pair of homological vector fields $Q_1$, $Q_2$ on $[-1]_{[-1]B}[-1]_A D \cong [-1]_{[-1]A}[-1]_B D$, such that the double Lie algebroid compatibility condition is equivalent to the condition $[Q_1,Q_2] = 0$.  It essentially follows from his result that the intermediate objects $[-1]_A D \to [-1]B$ and $[-1]_B D \to [-1]A$ are $Q$-algebroids.  Most of the examples of $Q$-algebroids that we will consider arise from double Lie algebroids in this manner.
\end{rmk}

\subsection{The $Q$-algebroid structure of $[-1]T\mathcal{A}$}\label{sec:oneta}

Let $\mathcal{A} \to \mathcal{M}$ be an algebroid.  Then
\begin{equation}\label{eqn:tadvb}
\xymatrix{T\mathcal{A} \ar[r] \ar[d] & \mathcal{A} \ar[d] \\ T\mathcal{M} \ar[r] & \mathcal{M}}
\end{equation}
is a double vector bundle in the sense of Pradines \cite{pradines}.  Applying the $[-1]$ functor to the rows results in the vector bundle $[-1]_\mathcal{A} T\mathcal{A} \to [-1]T\mathcal{M}$, where the subscript indicates the vector bundle structure with respect to which the degree-shift functor is taken.  It will now be shown that this vector bundle naturally has the structure of a $Q$-algebroid.  This is a supergeometric analogue of the fact in the ordinary case that (\ref{eqn:tadvb}) is a double Lie algebroid.

In what follows, when we write $[-1]T\mathcal{A}$, it will be assumed that the $[-1]$ functor is taken over $\mathcal{A}$.

The construction relies on the fact that the sections of $[-1]T\mathcal{A}$ are spanned by two types of ``lifts'' of sections of $\mathcal{A}$.  These lifts are analogous to the vertical and complete lifts of Yano and Ishihara \cite{yi}.

\begin{dfn} Let $X \in \Gamma(\mathcal{A})$.  The (odd) \emph{complete lift} $X^C$ of $X$ is a degree $|X|$ section of $[-1]T\mathcal{A}$ defined by the properties
\begin{align*} X^{C*} (\alpha) = X^* (\alpha), && X^{C*} (d\alpha) = d X^* (\alpha), \end{align*}
for $\alpha \in C^\infty(\mathcal{A})$.
\end{dfn}

\begin{dfn} Let $X \in \Gamma(\mathcal{A})$.  The (odd) \emph{vertical lift} $X^V$ of $X$ is a degree $|X|-1$ section of $[-1]T\mathcal{A}$ defined by the properties
\begin{align*} X^{V*} (\alpha) = 0, && X^{V*} (d\alpha) = (-1)^{|X| + |\alpha| + 1} X^* (\alpha), \end{align*}
for $\alpha \in C^\infty_{lin} (\mathcal{A})$.
\end{dfn}

\begin{prop} \label{prop:lifts}Let $f \in C^\infty(\mathcal{M})$, $X \in \Gamma(\mathcal{A})$.  Then
\begin{align*}
(fX)^C &= f \cdot X^C + (-1)^{|X| + |f|} df \cdot X^V, \\
(fX)^V &= f \cdot X^V.
\end{align*}
\end{prop}

\begin{proof}
The identities follow directly from the above definitions.
\end{proof}

\begin{prop} $\Gamma([-1]T\mathcal{A})$ is spanned by the complete and vertical lifts of $\Gamma(\mathcal{A})$.  More precisely, $\Gamma([-1]T\mathcal{A})$ is equal to $C^\infty([-1]T\mathcal{M}) \otimes \left( \left\{X^C\right\} \oplus \left\{X^V\right\}\right)$, modulo the relations of Proposition \ref{prop:lifts}.
\end{prop}

\begin{proof}Let $\{X_\alpha\}$ be a local frame of sections of $\mathcal{A}$ dual to fibre coordinates $\{\lambda^\alpha\}$.  Then it is clear from the definitions that $\{X^C_\alpha, -X^V_\alpha\}$ is a local frame of sections of $[-1]T\mathcal{A}$ dual to the fibre coordinates $\{\lambda^\alpha, \dot{\lambda}^\alpha\}$.
\end{proof}

The algebroid structure of $[-1]T\mathcal{A}$ is as follows.  The anchor $\widetilde{\rho} : [-1]T\mathcal{A} \to T([-1]T\mathcal{M})$ is defined by letting
\begin{align} \label{eqn:otanch}
\widetilde{\rho}(X^C) = \lie_{\rho(X)}, && \widetilde{\rho}(X^V) = \iota_{\rho(X)}, 
\end{align}
and then extending by $C^\infty([-1]T\mathcal{M})$-linearity.  Because the relations of Proposition \ref{prop:lifts} are identical to the relations satisfied by Lie derivative and contraction operators, $\widetilde{\rho}$ is well-defined.  The bracket is defined by letting
\begin{equation}
\begin{split} \label{eqn:otbrack}
\left[X^C, Y^C\right] &= \left[X,Y\right]^C, \\
\left[X^C, Y^V\right] &= \left[X,Y\right]^V, \\
\left[X^V, Y^V\right] &= 0,
\end{split}
\end{equation}
and extending by the Leibniz rule.  

One could now verify directly that the bracket (\ref{eqn:otbrack}) satisfies the Jacobi identity.  However, we will instead describe the associated operator $d_{[-1]T\mathcal{A}} \in \vect\left([-1]_{[-1]T\mathcal{M}}\left([-1]T\mathcal{A}\right)\right)$ and see that $\left(d_{[-1]T\mathcal{A}}\right)^2 = 0$.

As supermanifolds, and in fact as double vector bundles, $[-1]_{[-1]T\mathcal{M}}\left([-1]T\mathcal{A}\right)$ may be identified with $[-1]_{[-1]\mathcal{A}} T \left([-1]\mathcal{A}\right)$.  The algebroid structure on $\mathcal{A}$ is associated to a differential $d_\mathcal{A} \in \vect([-1]\mathcal{A})$ whose Lie derivative operator is $\lie_{d_\mathcal{A}} \in \vect\left([-1]_{[-1]\mathcal{A}}T ([-1]\mathcal{A})\right)$.

\begin{thm}\label{thm:lieda}$d_{[-1]T\mathcal{A}} = \lie_{d_\mathcal{A}}$.
\end{thm}
\begin{proof}
From the local coordinate description (\ref{eqn:dalocal}) of $d_\mathcal{A}$, the Lie derivative may be computed to be 
\begin{equation*} 
\begin{split}
\lie_{d_\mathcal{A}} = & \lambda^\alpha \rho_\alpha^i \pdiff{}{x^i} 
- (-1)^{p_\alpha (p_\beta - 1)}\frac{1}{2} \lambda^\alpha \lambda^\beta c_{\alpha \beta}^\gamma \pdiff{}{\lambda^\gamma} 
+ (-1)^{p_\alpha} \lambda^\alpha d\rho^i_\alpha \pdiff{}{\dot{x}^i} \\
&- (-1)^{p_\alpha p_\beta} \lambda^\alpha \dot{\lambda}^\beta c_{\alpha \beta}^{\gamma} \pdiff{}{\dot{\lambda}^\gamma}
- \dot{\lambda}^\alpha \rho_\alpha^i \pdiff{}{\dot{x}^i} 
+ (-1)^{(p_\alpha +1) p_\beta} \frac{1}{2} \lambda^\alpha \lambda^\beta dc_{\alpha \beta}^\gamma \pdiff{}{\dot{\lambda}^\gamma}.
\end{split}
\end{equation*}
If $\{X_\alpha\}$ are the sections dual to the fibre coordinates $\{\lambda^\alpha\}$ on $\mathcal{A}$, then $\{X^C_\alpha, - X^V_\alpha\}$ are the sections dual to $\{\lambda^\alpha,\dot{\lambda}^\alpha\}$ on $[-1]T\mathcal{A}$.  From equations (\ref{eqn:anchor}) and (\ref{eqn:bracket}), it is straightforward to check that (\ref{eqn:otanch}) and (\ref{eqn:otbrack}) are satisfied.
\end{proof}

\begin{cor}
\begin{enumerate}
\item $\left(d_{[-1]T\mathcal{A}}\right)^2 = 0$ and, equivalently, the bracket on $\Gamma([-1]T\mathcal{A})$ satisfies the Jacobi identity.
\item The de Rham differential $d \in \vect\left([-1]T([-1]\mathcal{A})\right)$ satisfies the equation $\left[d_{[-1]T\mathcal{A}}, d \right] = 0$.
\end{enumerate}
\end{cor}

\begin{cor}
$[-1]T\mathcal{A}$ is a $Q$-algebroid with the algebroid structure described above and the morphic vector field $d$.
\end{cor}

\subsection{Example: The Weil algebra}
\label{ex:weil}

Consider the simple example of a Lie algebra $\mathfrak{g}$.  The construction of \S\ref{sec:oneta} results in the structure of a \emph{Lie $Q$-algebra} (i.e.\ a $Q$-algebroid where the base is a point) on $[-1]T\mathfrak{g} = \mathfrak{g} \oplus [-1]\mathfrak{g}$.  To obtain the associated double complex, we apply the $[-1]$ functor to get $[-1]_{\{pt.\}}([-1]T\mathfrak{g}) = [-1]\mathfrak{g} \oplus [-2]\mathfrak{g}$.  The algebra $C^\infty([-1]\mathfrak{g} \oplus [-2]\mathfrak{g})$ is equal to the \emph{Weil algebra} $\mathcal{W}(\mathfrak{g}) \defequal \bigwedge \mathfrak{g}^* \otimes \bigS \mathfrak{g}^*$.

Let $\{v_i\}$ be a basis for $\mathfrak{g}$, and let $\{\theta^i\}$ and $\{\dot{\theta^i}\}$ be the dual bases in degree $1$ and $2$, respectively.  The de Rham differential on $[-1]T\mathfrak{g}$ induces the differential operator $d_K = \dot{\theta}^i \pdiff{}{\theta^i}$, which is known as the \emph{Koszul operator}.

In coordinates, The Lie algebra differential $d_{[-1]T\mathfrak{g}}$ is
\begin{equation*}
d_{[-1]T\mathfrak{g}} = -\frac{1}{2} c_{ij}^k \theta^i \theta^j \pdiff{}{\theta^k} - c_{ij}^k \theta^i \dot{\theta}^j \pdiff{}{\dot{\theta}^k},
\end{equation*}
and the total differential $d_{[-1]T\mathfrak{g}} + d_K$ is known as the Weil differential, denoted $d_\mathcal{W}$. 

The Weil algebra $\mathcal{W}(\mathfrak{g})$, equipped with the Weil differential, is known to be an acyclic complex.  However, if $G$ is a compact, connected Lie group with Lie algebra $\mathfrak{g}$, one can obtain a model for the cohomology of $BG$ by restricting to a certain subcomplex, defined as follows.

For any element $v = a^i v_i \in \mathfrak{g}$, there is a ``contraction operator'' $I_v = a^i \pdiff{}{\theta^i}$.

\begin{dfn} An element $\omega \in \mathcal{W}(\mathfrak{g})$ is called
\begin{enumerate}
\item \emph{horizontal} if $I_v \omega = 0$ for all $v \in \mathfrak{g}$,
\item \emph{invariant} if $L_v \omega \defequal [I_v, d_\mathcal{W}] \omega = 0$ for all $v \in \mathfrak{g}$, and
\item \emph{basic} if $\omega$ is both horizontal and invariant.
\end{enumerate}
\end{dfn}

It is clear that the horizontal subalgebra is $\bigS\mathfrak{g}^*$.  Furthermore, $d_\mathcal{W}$ vanishes on the basic subcomplex $\left(\bigS\mathfrak{g}^*\right)^G$.

The following theorem is a classic result due to Cartan.

\begin{thm}[\cite{cartan}]\label{thm:cartan}
If $G$ is a compact and connected Lie group, then cohomology of $BG$ is equal to the cohomology of the basic subcomplex of $\mathcal{W}(\mathfrak{g})$, which is equal to $\left(\bigS\mathfrak{g}^*\right)^G$.
\end{thm}

\subsection{The double complex of $[-1]TA$}\label{sec:doubleoneta}

We now consider the more general case of a $Q$-algebroid of the form $[-1]TA \to [-1]TM$, where $A \to M$ is a Lie algebroid.  Using the identification $[-1]_{[-1]TM}[-1]_A TA = [-1]_{[-1]A}T([-1]A)$, the space of cochains in the associated double complex can be identified with the space of differential forms on $[-1]A$.  It was shown in Theorem \ref{thm:lieda} that the algebroid differential of $[-1]TA$ is equal to $\lie_{d_A}$.  The total differential is the sum of the algebroid differential and the given morphic vector field, which in this case is the de Rham operator $d$.  Thus the total complex may be described simply as
\begin{equation}\label{eqn:totalcom}
\left(\Omega([-1]A), \lie_{d_A} + d \right).
\end{equation}

\begin{ex}[BRST model of equivariant cohomology, part I]\label{ex:brst1}

Let $M$ be a manifold and let $G$ be a compact, connected Lie group with a right action on $M$.  The infinitesimal data of the action may be described by the action algebroid $M \times \mathfrak{g} \to M$, and from this we obtain the $Q$-algebroid $[-1]T(M \times \mathfrak{g}) = [-1]TM \times [-1]T\mathfrak{g}$.  The algebra of cochains is $C^\infty\left([-1]_{[-1]TM}([-1]TM \times [-1]T\mathfrak{g})\right)$, which is naturally isomorphic to $\Omega(M) \otimes \mathcal{W}(\mathfrak{g})$.  The de Rham differential for $M \times \mathfrak{g}$ then splits into
$d_M + d_K$, where $d_M$ is the de Rham differential on $M$ and $d_K$ is the Koszul operator on $\mathcal{W}(\mathfrak{g})$.

The total differential $D_B$ may be written as
\begin{equation}\label{brstdiff}
D_B = d_M + d_{\mathcal{W}} + \theta^i \lie_{\rho(v_i)} - \dot{\theta}^i \iota_{\rho(v_i)},
\end{equation}
where $d_{\mathcal{W}}$ is the Weil differential and $\rho: \mathfrak{g} \to \vect(M)$ describes the infinitesimal action of $\mathfrak{g}$ on $M$.

The algebra $\Omega(M) \otimes \mathcal{W}(\mathfrak{g})$ and the differential (\ref{brstdiff}) form the \emph{BRST model} of equivariant cohomology \cite{kalkman:brst, osb}.  However, this complex is a model for $M \times EG$ (hence its cohomology is equal to $H^\bullet(M)$), and as in Example \ref{ex:weil} one must restrict to a suitably defined basic subcomplex in order to obtain the equivariant cohomology.  We will return to this issue in \S\ref{sec:basic}.
\end{ex}

\subsection{The Mathai-Quillen-Kalkman isomorphism}\label{sec:mathai-quillen}

Another model of equivariant cohomology, which is more well-known than the BRST model, is the Weil model.  The Weil model has the same algebra $\Omega(M) \otimes \mathcal{W}(\mathfrak{g})$ as the BRST model, but the simpler differential $d_M + d_K$.  In \cite{kalkman:brst}, Kalkman described an extension of the Mathai-Quillen isomorphism that relates the BRST differential and the Weil model differential, thus showing that the two models are equivalent\footnote{The original isomorphism of Mathai and Quillen \cite{mq} sent the Cartan model, which may be identified with the basic subcomplex of the BRST model, to the basic subcomplex of the Weil model.  Kalkman extended the isomorphism to the total complexes.}.  In order to compute the cohomology of the complex (\ref{eqn:totalcom}), we will use a generalization of the Mathai-Quillen-Kalkman isomorphism, described as follows.

Since $d_A$ is a degree $1$ vector field on $[-1]A$, it follows that $\iota_{d_A}$ is a degree $0$ vector field on $\Omega([-1]A)$; specifically, in terms of the double grading, $\iota_{d_A}$ is of degree $(1,-1)$.
\begin{dfn} The \emph{generalized Mathai-Quillen-Kalkman isomorphism} is $\gamma \defequal \exp (\iota_{d_A})$.
\end{dfn}

\begin{lemma} \label{lemma:mqkrel}
The de Rham differential $d$ is $\gamma$-related to $d + \lie_{d_A}$.
\end{lemma}
\begin{proof} Since $\iota_{d_A}$ is nilpotent, the identity $\Ad_{\exp (\iota_{d_A})} = \exp (\ad_{\iota_{d_A}})$ holds.  From the Cartan commutation relations, it is immediate that
\begin{align*}
\ad_{\iota_{d_A}} (d) &= [\iota_{d_A}, d] = \lie_{d_A}, \\
\ad^2_{\iota_{d_A}} (d) &= [\iota_{d_A}, \lie_{d_A}] = -\iota_{[d_A,d_A]} = 0,
\end{align*}
and it follows that $\Ad_{\gamma} (d) = d + \lie_{d_A}$.
\end{proof}

\begin{rmk} The differential for the Weil model is sometimes taken to be $d_M + d_\mathcal{W}$, where $d_\mathcal{W}$ is the Weil differential.  The isomorphism given by Kalkman in fact relates this differential to the BRST differential.  Thus, to compare his isomorphism with ours, it is necessary to use an automorphism of $\mathcal{W}(\mathfrak{g})$ that relates $d_K$ and $d_\mathcal{W}$.  Such an isomorphism is well-known (e.g.\ \cite{gs}), but we point out that it is a special case of $\gamma$, when $A = \mathfrak{g}$.
\end{rmk}

\begin{cor} \label{cor:eg} The $Q$-algebroid cohomology of $[-1]TA$ is equal to $H^\bullet(M)$.
\end{cor}

\begin{proof}
The isomorphism $\gamma$ provides an isomorphism of the total complex and the de Rham complex of $[-1]A$.  Using the Euler vector field for the vector bundle $[-1]A \to M$, it is a simple exercise to construct a chain homotopy between $\Omega([-1]A)$ and $\Omega(M)$.
\end{proof}

\subsection{The problem of the basic subcomplex}\label{sec:basic}
	
If $A$ is the Lie algebroid of a Lie groupoid $G \arrows M$, then the result of Corollary \ref{cor:eg} is consistent with the assertion that the complex (\ref{eqn:totalcom}) is a cohomological model for $EG$.  Following Examples \ref{ex:weil} and \ref{ex:brst1}, it seems reasonable to look for a basic subcomplex that computes the cohomology of $BG$.  In this section, we suggest a definition of such a basic subcomplex.  In the example of an action algebroid, this definition specializes to the existing definition of the basic subcomplex, which computes the equivariant cohomology.  However, the definition has a severe drawback, namely that it requires a choice of a connection on $A$.  The resulting basic cohomology is therefore vulnerable to the possibility of being dependent on the choice; for this reason, we choose, here and in \S\ref{sec:other}, to focus on examples where there is a canonical choice.

Recall (see (\ref{eqn:tadvb})) that $[-1]TA$ has the double vector bundle structure
\begin{equation}\label{eqn:onetadvb}
\dvb{[-1]TA}{A}{[-1]TM}{M}.
\end{equation}
Although $[-1]TA$ does not have a canonical vector bundle structure over $M$, it fits into the sequence
\begin{equation}\label{eqn:onetaexact}
[-1]A \to [-1]TA \to A \oplus [-1]TM,
\end{equation}
where $[-1]A$ is identified with the (degree-shifted) bundle of tangent vectors along the zero-section of $A$ that are tangent to the fibres.  In Mackenzie's \cite{mac:dblie1} terminology, $[-1]A$ is the \emph{core} of the double vector bundle (\ref{eqn:onetadvb}).  In fact, (\ref{eqn:onetaexact}) may be viewed as an exact sequence of double vector bundles
\begin{equation}\label{eqn:onetaexactdvb}
\dvb{[-1]A}{M}{\protect\vphantom{[-1]}M}{M} \raisebox{-4.5ex}{$\to$} \dvb{[-1]TA}{A}{[-1]TM}{M} \raisebox{-4.5ex}{$\to$} \dvb{A \oplus [-1]TM}{A}{[-1]TM}{M}.
\end{equation}

A section (in the category of double vector bundles) of the exact sequence (\ref{eqn:onetaexact}) is clearly equivalent to a splitting $[-1]TA \cong [-1]A \oplus A \oplus [-1]TM$.  The existence of such a section follows from

\begin{lemma}\label{lemma:splitting}
There is a one-to-one correspondence between sections of (\ref{eqn:onetaexact}) and linear connections on $A$.
\end{lemma}
\begin{proof}
If we disregard the $[-1]$'s, or more rigorously, eliminate them by applying the $[1]$ functor to the rows of (\ref{eqn:onetaexactdvb}), then we see that a section provides a way to lift vectors in $T_x M$, $x \in M$ to vectors in $T_{\tilde{x}}A$, where $\tilde{x}$ is a point in the fibre over $x$.  The requirement that the section be linear with respect to the horizontal vector bundle structures of (\ref{eqn:onetaexactdvb}) asserts that the lifts arise from linear maps $T_x M \to T_{\tilde{x}}A$, and the requirement that the section be linear with respect to the vertical structures assert that the lifts respect the linear structure of $A$.
\end{proof}

Let us now fix a connection, and therefore a splitting $[-1]TA \cong [-1]A \oplus A \oplus [-1]TM$.  In terms of the double complex of $[-1]TA$, the splitting gives an identification 
\begin{equation}\label{eqn:splitdblcx}
\Omega([-1]A) \cong \bigS \Gamma(A^*) \otimes \bigwedge \Gamma(A^*) \otimes \Omega(M)
\end{equation}
For any $X \in \Gamma(A)$, let $I_X$ denote the operator that acts by contraction in the exterior algebra component of (\ref{eqn:splitdblcx}).

\begin{dfn} An element $\omega \in C^\infty([-1]([-1]TA))$ is called
\begin{enumerate}
\item \emph{horizontal} if $I_X \omega = 0$ for all $X \in \Gamma(A)$,
\item \emph{invariant} if $L_X \omega \defequal [I_X, d + \lie_{d_A}] \omega = 0$ for all $X \in \Gamma(A)$, and
\item \emph{basic} if $\omega$ is both horizontal and invariant.
\end{enumerate}
\end{dfn}

\begin{rmk}\label{rmk:genconn}
The result of Lemma \ref{lemma:splitting} is due to Gracia-Saz and Mackenzie \cite{alfonso:dvb}.  In the general context of a double vector bundle (\ref{eqn:dvb}), they have shown that there always exist splittings $D \cong C \oplus A \oplus B$, where $C \to M$ is the core vector bundle.  They refers to such a splitting as a \emph{generalized connection}.  The definition of the basic subcomplex may be extended in the obvious way to any double complex arising from a double Lie algebroid equipped with a generalized connection.
\end{rmk}

\begin{ex}[BRST model, part II]\label{ex:brst2}
Recall the BRST model of Example \ref{ex:brst1}.  The action algebroid $M \times \mathfrak{g} \to M$ has a canonical flat connection and thus a natural splitting
\begin{equation*}
[-1]T(M \times [-1]\mathfrak{g}) = [-1]TM \times [-1]\mathfrak{g} \times [-2]\mathfrak{g},
\end{equation*}
or, equivalently,
\begin{equation*}
\Omega(M \times [-1]\mathfrak{g}) = \Omega(M) \otimes \bigwedge \mathfrak{g}^* \otimes \bigS \mathfrak{g}^*.
\end{equation*}

If $\{X_i\}$ is the global frame of flat sections corresponding to the basis $\{v_i\}$ of $\mathfrak{g}$,
then the contraction operators $I_X$ are generated freely as a $C^\infty(M)$-module by those of the form
\begin{equation*}
I_{X_i} = \pdiff{}{\theta^i}.
\end{equation*}
The horizontal elements are simply those that do not depend on any $\theta^i$, i.e.\ those that lie in $\Omega(M) \otimes S \mathfrak{g}^*$.  If we restrict $L_{x_i}$ to the horizontal subalgebra, we have
\begin{equation*}
L_{X_i} = \left[\pdiff{}{\theta^i}, D_B\right] = -\dot{\theta}^j c_{ij}^k \pdiff{}{\dot{\theta}^k} + \lie_{\rho(v_i)}.
\end{equation*}
Therefore the basic subalgebra is $\left(\Omega(M) \otimes S \mathfrak{g}^*\right)^G$, and on this subalgebra the differential becomes the Cartan differential
\begin{equation*}
d_C \defequal d_M - \dot{\theta}^i \iota_{\rho(v_i)}.
\end{equation*}
Thus the definition of basic elements agrees in this case with the usual one, where the basic subcomplex of the BRST model is equal to the Cartan model.
\end{ex}

\section{Other examples}\label{sec:other}

\subsection{Vacant double Lie algebroids}\label{sec:vacant}
	
A double Lie algebroid (\ref{eqn:dvb}) is said to be \emph{vacant} \cite{mac:dbl2} if the double-projection map $D \to A \oplus B$ is a diffeomorphism.  In other words, a vacant double Lie algebroid is of the form
\begin{equation}\label{eqn:dlavacant}
\dvb{A\oplus B}{A}{B}{M}.
\end{equation}
Applying the $[-1]$ functor to the rows, we obtain the $Q$-algebroid $A\oplus [-1]B \to [-1]B$.  Let $d_{\bar{A}} \in \vect([-1]A \oplus [-1]B)$ denote the algebroid differential.  By the properties of double Lie algebroids, we have that $d_{\bar{A}}$ is a linear vector field with respect to the bundle $[-1]A \oplus [-1]B \to [-1]A$, with base vector field $d_A$.  In coordinates $\{x^i, \alpha^i, \beta^i\}$, where $\{\alpha^i\}$ and $\{\beta^i\}$ are fibre coordinates on $[-1]A$ and $[-1]B$, respectively, $d_{\bar{A}}$ is therefore of the form
\begin{equation*}
d_{\bar{A}} = \alpha^i \left(\rho_i^j(x) \pdiff{}{x^j} + \sigma_{ij}^k(x) \beta^j \pdiff{}{\beta^k}\right) - c_{ij}^k(x) \alpha^i \alpha^j \pdiff{}{\alpha^k},
\end{equation*}
where $\rho_i^j$ and $c_{ij}^k$ are the anchor and structure functions for $A$.  The additional data $\sigma_{ij}^k$ describes a representation of $A$ on the vector bundle $B$.  

Similarly, the morphic vector field $d_{\bar{B}}$ is determined by the Lie algebroid structure of $B$ and a representation of $B$ on the vector bundle $A$.  The compatibility condition $[d_{\bar{A}}, d_{\bar{B}}] = 0$ is equivalent to the condition that $A$ and $B$, with their mutual representations, form a \emph{matched pair} in the sense of Mokri \cite{mokri}.  Furthermore, the total differential $d_{\bar{A}} + d_{\bar{B}}$ gives $[-1]A \oplus [-1]B$ the structure of an antialgebroid, corresponding to the Lie algebroid structure on $A \oplus B$ associated to the matched pair \cite{mokri}.

In this case, the space of cochains for the double complex is $\bigwedge \Gamma(A^*) \otimes \bigwedge \Gamma(B^*)$.  Because there is no core, no choice is required in order to define the basic subcomplex for a vacant double Lie algebroid.  The horizontal elements are those that vanish under contraction by sections of $A$; in other words, the horizontal subalgebra is $\bigwedge \Gamma(B^*)$.  The basic subcomplex consists of those elements of $\bigwedge \Gamma(B^*)$ that are invariant with respect to $A$.  On the basic subcomplex, the total differential coincides with $d_B$, so the basic cohomology is equal to the $A$-invariant Lie algebroid cohomology of $B$.  

We warn the reader that, although the notion of double Lie algebroid is symmetric with respect to the roles of $A$ and $B$, the definition of the basic subcomplex is not.  The reason is that the basic subcomplex is meant to be a model for the double complex \cite{m:q} of the $\LA$-groupoid
\begin{equation}\label{eqn:lavacant}
\la{s^*(B)}{G}{B}{M},
\end{equation}
where $G \arrows M$ is a Lie groupoid integrating $A$, and $(G,B)$ has a matched pair structure integrating that of $(A,B)$.  Indeed, if $G$ has compact, connected $t$-fibres, then the $\LA$-groupoid cohomology of (\ref{eqn:lavacant}) agrees with the basic cohomology of (\ref{eqn:dlavacant}).

Examples of vacant double Lie algebroids are discussed in \S\ref{sec:poissact} and \S\ref{sec:bialg} below.

\subsection{Poisson actions}\label{sec:poissact}
The following example was introduced by Lu \cite{lu} in relation to the study of Poisson homogeneous spaces (also see \cite{mac:reduction}).

Let $(M, \pi)$ be a Poisson manifold, and let $(\mathfrak{g}, \mathfrak{g}^*)$ be a Lie bialgebroid with a (right) Poisson\footnote{A right action of $\mathfrak{g}$ on $M$ is a Lie algebra morphism $\rho: \mathfrak{g} \to \vect(M)$.  The action is Poisson if the natural extension of $\rho$ to a map $\bigwedge \mathfrak{g} \to \vect^\bullet(M)$ is a morphism of differential Gerstenhaber algebras.} action $\rho: \mathfrak{g} \to \vect(M)$.  There is a corresponding representation of the action algebroid $M \times \mathfrak{g}$ on $T^*M$, defined on horizontal sections $a \in \mathfrak{g}$ by $a \mapsto \lie_{\rho(a)}$, where $\lie_{\rho(a)}$ is thought of as an operator on $\Gamma(T^*M)$.  

On the other hand, there is a representation of $T^*M$ on $M \times \mathfrak{g}$, defined as follows.  The anchor map $M \times \mathfrak{g} \to TM$, which, by abuse of notation, we will also refer to as $\rho$, dualizes to a map $\rho^*: \Omega(M) \to \Gamma(M \times \mathfrak{g}^*)$.  For $\alpha \in \Omega^1(M)$, $\rho^*(\alpha)$ acts on horizontal sections of $\Gamma(M \times \mathfrak{g})$ by the fibrewise coadjoint action of $\mathfrak{g}^*$ on $\mathfrak{g}$.  The action is extended to all sections by the property $\rho^*(\alpha)(fa) = \pi^\sharp \alpha (f) a + f \rho^*(\alpha)(a)$.

The mutual actions of $M \times \mathfrak{g}$ and $T^*M$ form a matched pair \cite{lu}, so there is a vacant double Lie algebroid of the form
\begin{equation}\label{eqn:poissact}
\dvb{(M \times \mathfrak{g}) \oplus T^*M}{M \times \mathfrak{g}}{T^*M}{M}.
\end{equation}
The associated double complex is $\bigwedge \mathfrak{g}^* \otimes \vect^\bullet(M)$, and the basic subcomplex computes the $\mathfrak{g}$-invariant Poisson cohomology of $M$.  This basic subcomplex is smaller than the one defined by Lu, which, in the case where the action is transitive, may be identified with the tensor product of the $\mathfrak{g}$-invariant de Rham complex and the $\mathfrak{g}$-invariant Lichnerowicz-Poisson complex.

The double Lie algebroid (\ref{eqn:poissact}) has also arisen in relation to Mackenzie's \cite{mac:reduction} description of Poisson reduction.  In this procedure, one passes to $K \subseteq T^*M$, where $K$ is the kernel of the bundle map $\mathfrak{p}: T^*M \to M \times \mathfrak{g}^*$, for which the associated map of sections is $\rho^*$.  If $\mathfrak{p}$ has constant rank, then $K$ is a subbundle, in which case it is also a Lie subalgebroid.  The action of $M \times \mathfrak{g}$ on $T^*M$ restricts to $K$, so one can form the vacant double Lie algebroid
\begin{equation*}
\dvb{(M \times \mathfrak{g}) \oplus K}{M \times \mathfrak{g}}{K}{M},
\end{equation*}
The basic subcomplex then consists of $\mathfrak{g}$-invariant elements of $\bigwedge \Gamma(K^*)$, which can be identified with $\left(\vect^\bullet(M)/ \left\langle \im \rho\right\rangle\right)^\mathfrak{g}$.  If $G$ is a connected Lie group with Lie algebra $\mathfrak{g}$ and if the quotient $M/G$ is a manifold, then the basic subcomplex is naturally isomorphic to the Lichnerowicz-Poisson complex of $M/G$.

\subsection{Lie bialgebras and Drinfel'd doubles}\label{sec:bialg}
Let $(\mathfrak{g}, \mathfrak{g}^*)$ be a Lie bialgebra.  One can form a double Lie algebroid
\begin{equation}\label{eqn:dlabialg}
\dvb{\mathfrak{g}\oplus \mathfrak{g}^*}{\mathfrak{g}^*}{\mathfrak{g}}{\cdot},
\end{equation}
where the two Lie algebroid structures of $\mathfrak{g}\oplus \mathfrak{g}^* = T^*\mathfrak{g} = T^*\mathfrak{g}^*$ correspond to the Poisson structures on $\mathfrak{g}$ and $\mathfrak{g}^*$, respectively.  Equivalently, they may be viewed as action algebroids for the coadjoint action of $\mathfrak{g}^*$ on $\mathfrak{g}$ and vice versa.

Applying the $[-1]$ functor to the rows, we obtain the $Q$-algebroid $\mathfrak{g}^* \oplus [-1]\mathfrak{g} \to [-1]\mathfrak{g}$, where the algebroid differential $\delta_* \in \vect\left([-1]\mathfrak{g}^* \oplus [-1]\mathfrak{g}\right)$ may be identified with the Chevalley-Eilenberg differential for $\mathfrak{g}^*$ with coefficients in $\bigwedge \mathfrak{g}^*$, and the morphic vector field $\delta \in \vect\left([-1]\mathfrak{g}^* \oplus [-1]\mathfrak{g}\right)$ may similarly be identified with the Chevalley-Eilenberg differential for $\mathfrak{g}$ with coefficients in $\bigwedge \mathfrak{g}$.  The property $\left[\delta, \delta_* \right] = 0$ is the homological version of the Lie bialgebra compatibility condition \cite{ks1, lr}.

The algebra of cochains for the double complex is $C^\infty([-1]\mathfrak{g} \oplus [-1]\mathfrak{g}^*) = \bigwedge \mathfrak{g}^* \otimes \bigwedge \mathfrak{g}$, and the total differential gives $[-1]\mathfrak{g} \oplus [-1]\mathfrak{g}^*$ a Lie antialgebra structure, for which the corresponding Lie algebra is the Drinfel'd \cite{drinfeld:ybe} double $\mathfrak{d} = \mathfrak{g} \oplus \mathfrak{g}^*$.

Since (\ref{eqn:dlabialg}) is a vacant double Lie algebroid, we have from \S\ref{sec:vacant} that the basic cohomology is the $\mathfrak{g}^*$-invariant Lie algebra cohomology of $\mathfrak{g}$.

\subsection{Lie bialgebroids and notions of double}

Much of the previous section generalizes to the case of Lie bialgebroids, in the following way.  If $(A, A^*)$ is a Lie bialgebroid, then one can form the double Lie algebroid
\begin{equation}\label{eqn:dlabialgbd}
\dvb{T^*A}{A^*}{A}{M}.
\end{equation}
This is the double Lie algebroid that Mackenzie \cite{mac:dbl} suggested as a generalization of the Drinfel'd double.  

Applying the $[-1]$ functor to the rows, we obtain the $Q$-algebroid $[-1]_{A^*} T^*A \to [-1]A$.  The algebroid differential is a homological vector field on $[-1]_A [-1]_{A^*} T^*A = [-2]_{[-1]A} T^*([-1]A)$, which is canonically symplectomorphic (via Roytenberg's \cite{royt} generalization of the Legendre tranform) to $[-2]_{[-1]A^*}T^*([-1]A^*)$.  With these identifications, we can identify the algebroid differential with the cotangent lift\footnote{The notion of cotangent lifts makes sense for degree-shifted cotangent bundles, in the following way.  A vector field $X \in \vect(\mathcal{M})$ may be viewed as a linear function on $[-n]T^*\mathcal{M}$, whose degree is $|X| + n$.  Since the Poisson bracket on $[-n]T^*\mathcal{M}$ is of degree $-n$, the Hamiltonian vector field $V_X$ is of degree $|X|$.} $V_{d_{A^*}}$ of $d_{A^*} \in \vect([-1]A^*)$, and we can identify the morphic vector field with the cotangent lift $V_{d_A}$ of $d_{A} \in \vect([-1]A)$.  The property $\left[V_{d_A}, V_{d_{A^*}} \right] = 0$ is the homological version of the Lie bialgebroid compatibility condition.  

In analogy with the case of Lie bialgebras, Roytenberg \cite{royt} proposed that $[-2]T^*([-1]A)$, equipped with the total differential $D \defequal V_{d_A} + V_{d_{A^*}}$, should be considered the Drinfel'd double of a Lie bialgebroid.  Thus we see that $Q$-algebroid cohomology provides a correspondence between Mackenzie's notion of double and that of Roytenberg\footnote{This correspondence was independently observed by Voronov \cite{voronov:mack}}.  
	
	\subsection{Equivariant Lie algebroid cohomology}\label{sec:equivalg}

As another example, we show that Ginzburg's \cite{ginzburg} notion of \emph{equivariant Lie algebroid cohomology} fits into the framework of $Q$-algebroids.  Let $A \to M$ be an algebroid and let $\mathfrak{g}$ be a Lie algebra. 

\begin{dfn} A (right) \emph{$A$-action} of $\mathfrak{g}$ on $M$ is a Lie algebra homomorphism $\tilde{a}: \mathfrak{g} \to \Gamma(A)$.
\end{dfn}

An $A$-action $\tilde{a}$ induces a Lie algebra homomorphism $a \defequal \rho \circ \tilde{a}: \mathfrak{g} \to \vect(M)$ that describes an action of $\mathfrak{g}$ on $M$.  If one begins with an action map $a: \mathfrak{g} \to \vect(M)$, then an $A$-action $\tilde{a}$ that lifts $a$ is an \emph{equivariant pre-momentum mapping} in the sense of Ginzburg \cite{ginzburg}.  An equivariant pre-momentum mapping, defined in this manner, is equivalent to an algebroid morphism from the action algebroid $M \times \mathfrak{g}$ to $A$.

\begin{thm}\label{thm:equivalg}
Let $\widetilde{a}: \mathfrak{g} \to \Gamma(A)$ be an $A$-action.  Then there is an induced Lie superalgebra action $\bar{\rho}: [-1]T\mathfrak{g} \to \vect([-1]A)$.  The associated action algebroid $[-1]A \times [-1]T\mathfrak{g} \to [-1]A$ is a $Q$-algebroid with morphic vector field $d_A + d_K$.
\end{thm}

\begin{proof}
Recall (\S\ref{sec:oneta}) that  $[-1]T\mathfrak{g}$ is naturally isomorphic to $\mathfrak{g} \oplus [-1]\mathfrak{g}$, where the two summands consist, respectively, of complete lifts $v^C$ and vertical lifts $v^V$ of elements $v \in \mathfrak{g}$.

The induced action is defined by
\begin{align}\label{eqn:barrho}
\bar{\rho}(v^V) = \iota_{\widetilde{a}(v)}, && \bar{\rho}(v^C) = \lie_{\widetilde{a}(v)} \defequal [\iota_{\widetilde{a}(v)}, d_A].
\end{align}
It follows from the Cartan relations that $\bar{\rho}$ is a Lie algebra homomorphism and thus describes an action of $[-1]T\mathfrak{g}$ on $[-1]A$.

Let $D$ be the degree $1$ operator on $\Gamma([-1]A \times [-1]T\mathfrak{g})$ defined on horizontal sections by
\begin{align} \label{eqn:onetgop} D(v^V) = v^C, && D(v^C) = 0 \end{align}
and extended to all sections by the property (see condition (1) of Lemma \ref{lem:morphderiv}) 
\begin{equation} \label{eqn:dextend} D(\omega X) = d_A \omega \cdot X + (-1)^{|\omega|} \omega DX,\end{equation}
 for $\omega \in C^\infty([-1]A)$ and $X \in [-1]T\mathfrak{g}$.  On horizontal sections, it is clear that $D$ is a derivation of the Lie bracket. In fact, (\ref{eqn:onetgop}) describes the operator on $[-1]T\mathfrak{g}$ corresponding to the Koszul vector field $d_K$.  In order to show that the extension of $D$ by (\ref{eqn:dextend}) is a derivation of the Lie bracket, it is sufficient to see that $\bar{\rho}(DX) = [d_A, \bar{\rho}(X)]$ for all $X \in [-1]T\mathfrak{g}$, which is immediate from (\ref{eqn:barrho}) and (\ref{eqn:onetgop}).  
 
By Theorem \ref{thm:morphderiv}, the operator $D$ corresponds to a morphic vector field $\Xi$ on $[-1]A \times [-1]T\mathfrak{g}$ with base vector field $d_A$.  Since, under the projection map $[-1]A \times [-1]T\mathfrak{g} \to [-1]T\mathfrak{g}$, $\Xi$ is related to $d_K$, we have that $\Xi = d_A + d_K$.
\end{proof}

Theorem \ref{thm:equivalg} gives us a double complex structure on the algebra $\bigwedge \Gamma(A^*) \otimes \mathcal{W}(\mathfrak{g})$.  As in \S\ref{sec:mathai-quillen}, we use a Mathai-Quillen-Kalkman-type isomorphism in order to compute the cohomology of the double complex.  

Let $d_K^*$ denote the vector field on $[-1]T\mathfrak{g}$ which, in the coordinates of \S\ref{ex:weil}, is given by $d_K^* = \theta^i \pdiff{}{\dot{\theta^i}}$.  The key property that $d_K^*$ satisfies is $[d_K^*, d_K] = e$, where $e$ is the Euler vector field on the vector space $[-1]T\mathfrak{g}$.  Let $d_{\bar{A}} \defequal d_{[-1]A \times [-1]T\mathfrak{g}}$ denote the differential for the algebroid $[-1]A \times [-1]T\mathfrak{g} \to [-1]A$.  We define $Q \defequal [d_K^*, d_{\bar{A}}]$.

\begin{lemma}$d_A + d_K$ is $\exp(-Q)$-related to the total differential $d_A + d_K + d_{\bar{A}}$.
\end{lemma}
\begin{proof}
As in Lemma \ref{lemma:mqkrel}, the nilpotency of $Q$ implies that $\Ad_{\exp(-Q)} = \exp(\ad_{-Q})$.  Thus we compute
\begin{equation*}\begin{split}
\ad_{-Q} (d_A + d_K) &= - \left[[d_K^*, d_{\bar{A}}], d_A + d_K\right] \\
&= -\left[d_K^*, [d_{\bar{A}}, d_A + d_K]\right] + \left[[d_K^*, d_A + d_K], d_{\bar{A}}\right].
\end{split}
\end{equation*}
The first term vanishes since $d_A + d_K$ is a morphic vector field.  Since $d_K^*$ commutes with $d_A$, the remaining term becomes $[e,d_{\bar{A}}] = d_{\bar{A}}$.  We then see that $\ad_{-Q}^2 (d_A + d_K) = 0$ since $d_{\bar{A}}^2 = 0$.  It follows that $\Ad_{\exp(-Q)}(d_A + d_K) = d_A + d_K + d_{\bar{A}}$.
\end{proof}

Because $\mathcal{W}(\mathfrak{g})$ is acyclic with respect to $d_K$, we have
\begin{cor}
The cohomology of the complex $\left(\bigwedge \Gamma(A^*) \otimes \mathcal{W}(\mathfrak{g}), d_A + d_K + d_{\bar{A}}\right)$ is equal to the algebroid cohomology of $A$.
\end{cor}

As usual, it is necessary to pass to a basic subcomplex in order to get more interesting cohomology.  In this case, there is a canonical splitting $[-1]A \times [-1]T\mathfrak{g} \cong [-1]A \times \mathfrak{g} \times [-1]\mathfrak{g}$, which is a generalized connection (see Remark \ref{rmk:genconn}) for the double vector bundle
\begin{equation*}
\dvb{[-1]A \times [-1]T\mathfrak{g}}{M \times \mathfrak{g}}{[-1]A}{M}.
\end{equation*}
Thus for any $v \in \mathfrak{g}$, there is a well-defined operator $I_v$, acting on $\bigwedge \Gamma(A^*) \otimes \bigwedge \mathfrak{g}^* \otimes \bigS \mathfrak{g}^*$ by contraction in the exterior algebra component.  The basic elements are those that are annihilated by $I_v$ and $L_v \defequal [I_v, d_A + d_K + d_{\bar{A}}]$ for all $v \in \mathfrak{g}$.  Generalizing the results of Example \ref{ex:brst2}, we have that the basic subalgebra is $\left(\bigwedge \Gamma(A^*) \otimes S(\mathfrak{g}^*) \right)^\mathfrak{g}$, and on this subalgebra the total differential becomes
\begin{equation}\label{eqn:pcdiff}
d_C = d_A - \dot{\theta}^i \iota_{\widetilde{a}(v_i)}.
\end{equation}
The basic subcomplex equipped with the differential (\ref{eqn:pcdiff}) is identical Ginzburg's \cite{ginzburg} model for equivariant algebroid cohomology.

\subsection{Lifted actions on Courant algebroids}

Let $E \to M$ be a Courant algebroid, where the inner product and (non-skew-symmetric) bracket are denoted by $\left\langle \cdot,\cdot\right\rangle$ and $\llbracket \cdot, \cdot \rrbracket$, respectively.  We will briefly sketch how a ``lifted action'', in the sense of Bursztyn, Cavalcanti, and Gualtieri \cite{bcg}, of a Lie algebra $\mathfrak{g}$ on $E$ induces a $Q$-algebroid structure on a corresponding action algebroid.  This provides an example of a $Q$-algebroid that does not originate from a double Lie algebroid.

Recall that Roytenberg \cite{royt:graded} has shown that there is a one-to-one correspondence between Courant algebroids and degree $2$ symplectic $Q$-manifolds.  The $Q$-manifold $\mathcal{E}$ associated to $E$ is (noncanonically) of the form $[-1]E^* \oplus [-2]T^*M$.  In particular, the space of degree $1$ functions on $\mathcal{E}$ is (canonically) equal to $\Gamma(E)$.  

Let $\mathfrak{g}$ be a Lie algebra with an action $a: \mathfrak{g} \to \vect(M)$.  A \emph{lifted action}\cite{bcg} is a map $\tilde{a}: \mathfrak{g} \to \Gamma(E)$ that respects the brackets and projects, via the anchor map, to the action map $a$.

\begin{thm}
Let $\widetilde{a}: \mathfrak{g} \to \Gamma(E)$ be a lifted action with isotropic image.  Then there is an induced Lie superalgebra action $\bar{\rho}: [-1]T\mathfrak{g} \to \vect(\mathcal{E})$.  The associated action algebroid $\mathcal{E} \times [-1]T\mathfrak{g} \to \mathcal{E}$ is a $Q$-algebroid.
\end{thm}

\begin{proof}
On $\{\mbox{degree $1$ functions}\} = \Gamma(E)$, the action is defined as follows:
\begin{align*}
\bar{\rho}(v^V) &= \left\langle \tilde{a}(v), \cdot \right\rangle \\
\bar{\rho}(v^C) &= \llbracket \tilde{a}(v), \cdot \rrbracket \defequal \ad_{\tilde{a}(v)}.
\end{align*}
More precisely, $\bar{\rho}(v^V)$ and $\bar{\rho}(v^C)$ are the Hamiltonian vector fields of $\tilde{a}(v)$ and $\delta \tilde{a}(v)$, respectively, where $\delta$ is the homological vector field on $\mathcal{E}$.  This can be shown to be a Lie algebra homomorphism (we note, in particular, that the hypothesis that $\im \tilde{a}$ be isotropic is used to show that $[\bar{\rho}(v^V), \bar{\rho}(w^V)] = 0$ for all $v,w \in \mathfrak{g}$).

As in the proof of Theorem \ref{thm:equivalg}, we define the operator $D$ on $\Gamma(\mathcal{E} \times [-1]T\mathfrak{g})$ by (\ref{eqn:onetgop}) and the property
\begin{equation} D(f X) = \delta f \cdot X + (-1)^{|f|} \omega DX,\end{equation}
for $f \in C^\infty(\mathcal{E})$ and $X \in [-1]T\mathfrak{g}$.  To see that $D$ is a derivation of the Lie bracket, it is sufficient to check that $\bar{\rho}(DX) = [\delta, \bar{\rho}(X)]$, which follows from the fact that $\delta$ is a derivation of the Poisson bracket and homological.  It follows that $D$ corresponds to a morphic vector field, which in this case is $\delta + d_K$.
\end{proof}

The reader may note the similarity between this construction and that of \S\ref{sec:equivalg}.  Following that case, we may define a basic subcomplex of $C^\infty(\mathcal{E}) \otimes \mathcal{W}(\mathfrak{g})$, the cohomology of which might reasonably be called \emph{equivariant Courant algebroid cohomology}.

\bibliographystyle{abbrv}
\bibliography{bibio}
\vspace{0.5cm}
\noindent
\address{\small\textsc{Rajan Amit Mehta\\
Instituto de Matem\'{a}tica Pura e Aplicada\\
Estrada Dona Castorina 110\\
Rio de Janeiro 22460-320\\
Brasil\\}}

\end{document}